\documentstyle{amltd}
\begin{document}
\annalsline{157}{2003}
\received{August 16, 2001}
\startingpage{647}
\def\bye{\end{document}}
 \font\tenrm=cmr10
\def\ritem#1{\item[{\rm #1}]}
\def\dfrac#1#2{\displaystyle\frac{#1}{#2}}
\catcode`\@=11
\font\twelvemsb=msbm10 scaled 1100
\font\tenmsb=msbm10
\font\ninemsb=msbm10 scaled 800
\newfam\msbfam
\textfont\msbfam=\twelvemsb  \scriptfont\msbfam=\ninemsb
  \scriptscriptfont\msbfam=\ninemsb
\def\msb@{\hexnumber@\msbfam}
\def\Bbb{\relax\ifmmode\let\next\Bbb@\else
 \def\next{\errmessage{Use \string\Bbb\space only in math
mode}}\fi\next}
\def\Bbb@#1{{\Bbb@@{#1}}}
\def\Bbb@@#1{\fam\msbfam#1}
\catcode`\@=12

 \catcode`\@=11
\font\twelveeuf=eufm10 scaled 1100
\font\teneuf=eufm10
\font\nineeuf=eufm7 scaled 1100
\newfam\euffam
\textfont\euffam=\twelveeuf  \scriptfont\euffam=\teneuf
  \scriptscriptfont\euffam=\nineeuf
\def\euf@{\hexnumber@\euffam}
\def\frak{\relax\ifmmode\let\next\frak@\else
 \def\next{\errmessage{Use \string\frak\space only in math
mode}}\fi\next}
\def\frak@#1{{\frak@@{#1}}}
\def\frak@@#1{\fam\euffam#1}
\catcode`\@=12


\newcommand{\thmref}[1]{Theorem~\ref{#1}}
\newcommand{\secref}[1]{\S\ref{#1}}
\newcommand{\lemref}[1]{Lemma~\ref{#1}}

\title{The best constant for the centered\\ Hardy-Littlewood maximal
inequality} 
\shorttitle{The best constant in maximal inequality} 

\author{Antonios D. Melas}
 
\institutions{University of Athens,  Athens, Greece\\
{\eightpoint {\it E-mail address\/}:
amelas@math.uoa.gr}}
 
 \centerline{\bf Abstract}
\vglue12pt 
We find the exact value of the best possible constant $C$ for the weak-type
$(1,1)$ inequality for the one-dimensional centered Hardy-Littlewood maximal
operator. We prove that $C$ is the largest root of the quadratic equation
$12C^{2}-22C+5=0$ thus obtaining $C=1.5675208\ldots \  $. This is the first time the
best constant for one of the fundamental inequalities satisfied by a
{\it centered} maximal operator is precisely evaluated.

\section{Introduction}

Maximal operators play a central role in the theory of differentiation of
functions and also in Complex and Harmonic Analysis. In general one considers
a certain collection of sets ${\cal C}$ in ${\Bbb  R}^{n}$ and then given
any locally integrable function~$f$, at each $x$ one measures the maximal
average value of $f$ with respect to the collection ${\cal C}$, translated
by $x$. Then it is of fundamental importance to obtain certain regularity
properties of this operators such as weak-type inequalities or $L^{p}$-boundedness. These properties are well known
if
${\cal C}$, for example, consists of all $\alpha D$ where $\alpha>0$ is arbritrary and $D\subseteq
{\Bbb  R}^{n}$ is a fixed bounded convex set containing $0$ in its interior.
Such maximal operators are usually called {\it centered}.

However little is known about the deeper properties of centered maximal
operators even in the simplest cases. And one way to acquire such a deeper
understanding is to start asking for the best constants in the corresponding
inequalities satisfied by them. In this direction let us mention the result of\break 
E. M. Stein and J.-O. Str\"{o}mberg \cite{Stein} where certain upper bounds are
given for such constants in the case of centered maximal operators as
described above, and the corresponding still open question raised there (see
also \cite[Problem 7.74b]{Hay}), on whether the best constant in the weak-type
$(1,1)$ inequality for certain centered maximal operators in ${\Bbb  R}^{n}$
has an upper bound independent of~$n$.

The simplest example of such a maximal operator is the centered
Hardy-Littlewood maximal operator defined by
\begin{equation}
M{}\,f{}(x)=\sup_{h>0}\frac{1}{2h}\int_{x-h}^{x+h}\left|  f\right| \label{i1}%
\end{equation}
for every $f\in L^{1}({\Bbb  R})$. The weak-type $(1,1)$ inequality for this
operator says that there exists a constant $C>0$ such that for every $f\in
L^{1}({\Bbb  R})$ and every $\lambda>0,$%
\begin{equation}
\left|  \{M\,f>\lambda\}\right|  \leq\frac{C}{\lambda}\left\|  f\right\|
_{1}.\label{i2}%
\end{equation}
However even in this case not much was known for the best constant
$C$ in the above inequality. This must be contrasted with the corresponding
uncentered maximal operator defined similarly to (\ref{i1}) but by not
requiring $x$ to be the center but just any point of the interval of
integration. Here the best constant in the analogous to (\ref{i2}) inequality
is equal to $2$ which corresponds to a single dirac delta. The proof follows
from a covering lemma that depends on a simple topological property of the
intervals of the real line and can be extended to the case of any measure of
integration, not just the Lebesgue measure (see \cite{Be}). Moreover in this
case the best constants in the corresponding $L^{p}$ inequalities are also
known (see \cite{Gra}).

However in the case of the centered maximal operator the behavior is much more
difficult and it seems to not only depend on the Lebesgue measure but to also
involve a much deeper geometry of the real line. A. Carbery proposed that
$C=3/2$ (\cite[Problem 7.74c]{Hay}), a joint conjecture with F. Soria  which
also appears in \cite{Tri} and corresponds to sums of equidistributed dirac
deltas. This conjecture has been refuted by J. M. Aldaz in \cite{Al} who
actually obtained the bounds $1.541\ldots =\dfrac{37}{24}\leq C\leq\dfrac
{9+\sqrt{41}}{8}=1.9253905\ldots <2$ which also implies that $C$ is strictly less
than the constant in the uncentered case, thus answering a question that was
asked in \cite{Tri}. Then J. Manfredi and F. Soria improved the lower bound
proving that (\cite{Man}; see also \cite{Al}): $C\geq\dfrac{5}{3}-\dfrac
{2\sqrt{7}}{3}\sin\left(\dfrac{\arctan(3\sqrt{3})^{-1}}{3}\right)=1.5549581\ldots \
$.
\vglue3pt
The proofs of these results use as a starting point the discretization
technique introduced by M. de Guzm\'{a}n \cite{Gu} as sharpened by M. Trinidad
Men\'{a}rguez-F. Soria (see Theorem 1 in \cite{Tri}). To describe it we
define for any finite measure $\sigma$ on ${\Bbb  R}$ the corresponding
maximal function
\begin{equation}
M{}\,\sigma{}(x)=\sup_{h>0}\frac{1}{2h}\int_{x-h}^{x+h}\left|  d\sigma\right|
.%
\end{equation}

Then the best constant $C$ in inequality (\ref{i2}) is equal to the
corresponding best constant in the inequality
\begin{equation}
\left|  \{M\,\mu>\lambda\}\right|  \leq\frac{C}{\lambda}\int_{{\Bbb  R}}%
d\mu\label{i4}%
\end{equation}

\noindent where $\lambda>0$ and $\mu$ runs through all measures of the form
$\sum_{i=1}^{n}\delta_{t_{i}}$ where $n\geq1$ and $t_{1},\ldots ,t_{n}%
\in{\Bbb  R}$. This technique allows us to apply arguments of combinatorial
nature to get information or bounds for this constant.

The author (see \cite{Mel}) using also this technique, obtained the following
improved estimates for $C$:
\begin{equation}
1.5675208\ldots =\frac{11+\sqrt{61}}{12}\leq C\leq\frac{5}{3}=1.66\ldots \label{i2a}%
\end{equation}
and also made the conjecture that the lower bound in (\ref{i2a}) is actually
the exact value of $C$. Recently in \cite{Mel1} the author found the best
constant in a related but more general covering problem on the real line. This
implies the following improvement of the upper bound in (\ref{i2a}):
$C\leq1+\dfrac{1}{\sqrt{3}}=1.57735\ldots \ $. None of these however tells us what
the exact value of $C$ is.

In this paper we will prove that the above conjecture is correct thus settling
the problem of the computation of the best constant $C$ completely. We will
prove the following.

\specialnumber{1}\proclaim{Theorem}
For the centered Hardy\/{\rm -}\/Littlewood maximal operator $M${\rm ,} for every measure
$\mu$ of the form $k_{1}\delta_{y_{1}}+\cdots +k_{y_{n}}\delta_{y_{n}}$ where
$k_{i}>0$ for $i=1,\ldots,n$ and $y_{1}<\cdots<y_{n}$ and for every $\lambda>0$ we
have
\begin{equation}
\left|  \{M\,\mu>\lambda\}\right|  \leq\frac{11+\sqrt{61}}{12\lambda}\left\|
\mu\right\| \label{mi}%
\end{equation}
and this is sharp.
\endproclaim

We will call the measures $\mu$ that appear in the statement of the above
theorem, {\it positive linear combinations of dirac deltas}.

In view of the discretization technique described above Theorem 1 implies the following.

\specialnumber{1}\proclaim{{C}orollary}
For every $f\in L^{1}({\Bbb  R})$ and for every $\lambda>0$ we have
\begin{equation}
\left|  \{M\,f>\lambda\}\right|  \leq\frac{11+\sqrt{61}}{12\lambda}\left\|
f\right\|  _{1}\label{mi1}%
\end{equation}
and this is sharp.
\endproclaim

Hence
\begin{equation}
C=\dfrac{11+\sqrt{61}}{12}=1.5675208\ldots \label{i4a}%
\end{equation}
is the largest solution of the quadratic equation
\begin{equation}
12C^{2}-22C+5=0.\label{i4b}%
\end{equation}
By the lower bound in (\ref{i2a}) proved in \cite{Mel} we only have to prove
inequality (\ref{mi}) to complete the proof of Theorem 1. The number appearing
in equality (\ref{i4a}) is probably not suggesting anything, nor is the
equation (\ref{i4b}). However this number is what one would get in the limit
by computing the corresponding constants in the measures that are produced by
applying an iteration based on the construction in \cite{Mel} that leads to
the lower bound. These measures, although rather complicated (much more
complicated than single or equidistributed dirac deltas), have a very distinct
inherent structure (see the appendix here). Thus it would be probably better
to view Theorem 1 as a statement saying that this specific structure actually
is one that produces configurations with optimal behavior.

Then, in a completely analogous manner as the result in \cite{Gu}, \cite{Tri},
we will also prove the following.

\specialnumber{2}\proclaim{Theorem}
For any finite Borel measure $\sigma$ on ${\Bbb  R}$ and for any $\lambda>0$
we have
\begin{equation}
\left|  \{M\,\sigma>\lambda\}\right|  \leq\frac{11+\sqrt{61}}{12\lambda
}\left\|  \sigma\right\|  .\label{mi2}%
\end{equation}
\endproclaim

We have included this here because it is then natural to ask whether there
exists a function $f\in L^{1}({\Bbb  R})$, or more generally a measure
$\sigma$, and a $\lambda>0$ for which equality holds in the corresponding
estimate (\ref{mi1}) and (\ref{mi2}). We will show here that such an extremal
cannot be found in the class of all positive linear combinations of dirac deltas.

\specialnumber{3}\proclaim{Theorem}
For any measure $\mu$ that is a positive linear combination of dirac deltas
and for any $\lambda>0$ we have
\begin{equation}
\left|  \{M\,\mu>\lambda\}\right|  <\frac{11+\sqrt{61}}{12\lambda}\left\|
\mu\right\|  .\label{mi3}%
\end{equation}
\endproclaim

For the proof of Theorem 1, that is of inequality (\ref{mi}), our starting
point will be the related covering and overlapping problems that were
introduced in \cite{Mel} using the discretization technique. This proof is
divided into several sections and will contain a mixture of combinatorial,
geometric and analytic arguments. We start from the assumption that this upper
bound is not correct and fix a certain combination of dirac deltas that
violates it and contain the least possible number of positions. Then using the
related covering problem from \cite{Mel}, studied in more detail here, we will
prove that this assumed measure will contain, or can be used to produce,
segments that share certain structural similarities with the examples leading
to the lower bound. This needs some work and is better described if we further
discretize the corresponding covering problem by assuming that all masses and
positions of this measure are integers. Then elaborating on the structure of
these segments combined with the assumed violation of (\ref{mi}) we will
obtain a certain estimate for the central part of these segments. This
estimate will then lead to a contradiction using the assumption that any
measure of fewer positions will actually satisfy (\ref{mi}). This will
complete the proof of Theorem 1. Then we will give the proofs of Theorems 2
and 3 and in the Appendix we will briefly describe the construction from
\cite{Mel} that leads to the lower bound and we will compare it with the proof
of the upper bound.

\demo{Acknowledgements}\hskip-4pt The author would like to thank Professors   A.
Carbery, L. Grafakos, J.-P. Kahane and F. Soria for their interest in this work.
\enddemo

\section{ Preliminaries}

We will start here by describing our basic reduction of the problem as was
introduced in \cite{Mel}, where also further details and proofs can be found.
We will consider measures $\mu$ of the form
\begin{equation}
\mu=\sum_{i=1}^{n}k_{i}\delta_{y_{i}}\label{a1}%
\end{equation}
where $n$ is a positive integer, $k_{1},\ldots ,k_{n}>0$ are its
{\it masses }and $y_{1}<\cdots <y_{n}$ are its {\it positions}.

For any such measure as in (\ref{a1}) we define the intervals
\begin{equation}
I_{i,j}=I_{i,j}(\mu)=\,[\,y_{j}-k_{i}-\cdots-k_{j}\,,\,y_{i}+k_{i}+\cdots+k_{j}%
\,],%
\end{equation}
for $1\leq i\leq j\leq n$ (where $[a,b]=\emptyset$ if $b<a$) and the set
\begin{equation}
E\,(\mu)=\bigcup_{1\leq i\leq j\leq n}I_{ij}(\mu)\,.%
\end{equation}
This set can be seen to be equal to $\{x:\,M\,\mu(x)\geq1/2\}$ (see
\cite{Mel}\noindent).

It will be convenient throughout this paper to use the following notation: We
define
\begin{equation}
K_{i}^{j}=k_{i}+\cdots+k_{j}\label{a1a}%
\end{equation}
if $1\leq i<j\leq n$, $K_{i}^{i}=k_{i}$ if $1\leq i\leq n$ and $K_{i}^{j}=0$
if $j<i$. Thus we can write $I_{i,j}(\mu)=\,[\,y_{j}-K_{i}^{j},\,y_{i}%
+K_{i}^{j}]$.

We will say that $\mu$ satisfies the {\it separability inequalities} if:
\begin{equation}
y_{i+1}-y_{i}>k_{i}+k_{i+1}\label{a2}%
\end{equation}
for all $i=1,\ldots ,n-1$. If this happens then it is easy to see that for any
$1\leq i<j\leq n$ we have
\begin{equation}
I_{i,j}(\mu)\subseteq(y_{i},y_{j}) \label{a3}%
\end{equation}
(in fact this is equivalent to $K_{i}^{j}<y_{j}-y_{i}$ which follows by adding
certain inequalities from (\ref{a2})) and therefore $E(\mu)\subseteq\lbrack
y_{1}-k_{1},y_{n}+k_{n}]$.

We also set
\begin{equation}
R\,(\mu)=\dfrac{\left|  E(\mu)\right|  }{2\left\|  \mu\right\|  }%
=\frac{\left|  E(\mu)\right|  }{2(k_{1}+\cdots + k_{n})}=\frac{\left|
E(\mu)\right|  }{2K_{1}^{n}}.%
\end{equation}
\noindent Then we have the following (see \cite{Mel}).

\specialnumber{1}\proclaim{Proposition}
{\rm (i)} The best constant $C$ in the Hardy\/{\rm -}\/Littlewood maximal inequality
{\rm (\ref{i2})} is equal to the supremum of all numbers $R\,(\mu)$ when $\mu$ runs
through all positive measures of the form {\rm (\ref{a1})} that satisfy {\rm (\ref{a2}).}
\vglue4pt
{\rm (ii)} $C$ is also equal to the supremum of all numbers $R\,(\mu)$ when $\mu$
runs through all positive measures as in {\rm (i)} that also satisfy the condition\/{\rm :}\/
\begin{equation}
E(\mu)=[y_{1}-k_{1},y_{n}+k_{n}].%
\end{equation}
\endproclaim

Any such measure that satisfies the conditions in Proposition 1(ii), that is
the separability inequalities and the connectedness of $E(\mu)$, will be
called {\it admissible}. It is clear that for any admissible $\mu$ the
intervals $I_{i,j}(\mu)$, $1\leq i\leq j\leq n$ form a covering of the
interval $[y_{1}-k_{1},y_{n}+k_{n}]$.

We will also use the following lemma whose proof is essentialy given in
\cite{Mel} (see also \cite{Al}).

\specialnumber{1}\proclaim{Lemma}
Suppose $\mu$ is a measure containing $n\geq2$ positions that does not satisfy
all separability inequalities {\rm (\ref{a2}),} that is for at least one $i $ we
have $y_{i+1}-y_{i}\leq k_{i+1}+k_{i}$. Then there exists an admissible
measure $\mu^{\ast}$ containing at most $n-1$ positions and such that
$R(\mu^{\ast})\geq R(\mu)$.
\endproclaim

Hence, unless otherwise stated, we will only consider measures $\mu$ that
satisfy all inequalities (\ref{a2}). It is easy then to see that for any such
$\mu$ the intervals $I_{i,i}(\mu)$ for $1\leq i\leq n$ are pairwise disjoint.
We define the set of {\it covered} {\it gaps }of $\mu$ as follows:
\begin{equation}
G(\mu)=E(\mu)\backslash\bigcup_{i=1}^{n}I_{i,i}(\mu).%
\end{equation}
This is the set of points that must be covered by the intervals $I_{i,j}(\mu)$
for $i<j$ that come from interactions of distant masses and are nonempty if
their positions are, in some sense, close together. We also have
\begin{equation}
R(\mu)=1+\dfrac{\left|  G(\mu)\right|  }{2K_{1}^{n}}.%
\end{equation}

To proceed further let us now fix an admissible measure $\mu$ as in
(\ref{a1}). An important device that can describe efficiently the covering
properties $I_{i,j}(\mu)$ for $i<j$ is the so called {\it gap interval }of
$\mu$ that was introduced in \cite{Mel}. We consider the positive numbers
\begin{equation}
x_{i}=y_{i+1}-y_{i}-k_{i+1}-k_{i}\label{b1}%
\end{equation}
for $1\leq i\leq n$, the points
\begin{equation}
a_{1}=0,\,a_{2}=x_{1},\,a_{3}=x_{1}+x_{2},\ldots , a_{n}=x_{1}+\cdots +
x_{n-1}\hskip.5in
\label{b2}%
\end{equation}
and define the gap interval $J(\mu)$ of $\mu$ as follows
\begin{equation}
J(\mu)=[a_{1},a_{n}].\label{b3}%
\end{equation}
The gap interval can be obtained from $E(\mu)=[y_{1}-k_{1},y_{n}-k_{n}]$ by
collapsing the central intervals $I_{i,i}(\mu)=[y_{i}-k_{i},y_{i}+k_{i}%
]$,\ $1\leq i\leq n$ into the points $a_{i}$. This can be described by
defining a (measure-preserving and discontinuous) mapping
\begin{equation}
Q=Q_{\mu}:J(\mu)\rightarrow G(\mu)\label{b3a}%
\end{equation}
that satisfies $Q(x)=y_{i}+k_{i}+(x-a_{i})$ whenever $x\in(a_{i},a_{i+1})$,
$1\leq i<n$. Thus $Q$ maps each subinterval $(a_{i},a_{i+1})$ of $J(\mu)$ onto
the corresponding gap $(y_{i}+k_{i},y_{i+1}-k_{i+1})$ of $G(\mu)$. It is also
trivial to see that the mapping $Q$ is distance nondecreasing and so $Q^{-1}$
is distance nonincreasing.

We also consider the intervals
\begin{equation}
J_{i}=J_{i}(\mu)=[a_{i}-k_{i},a_{i}+k_{i}]
\end{equation}
around each of the points $a_{i}$,$\ 1\leq i\leq n$, of $J(\mu)$, let
\begin{equation}
{\cal F}(\mu)=\{J_{1}(\mu),\ldots , J_{n}(\mu)\}
\end{equation}
denote the corresponding family of all these intervals and let
\begin{equation}
J_{i}^{+}=J_{i}^{+}(\mu)=[a_{i},a_{i}+k_{i}]\hbox{ and }J_{i}^{-}=J_{i}%
^{-}(\mu)=[a_{i}-k_{i},a_{i}]\hskip.5in
\end{equation}
denote the right and left half of $J_{i}$ respectively. We also consider the
families of intervals
\begin{equation}
{\cal F}^{+}(\mu)=\{J_{1}^{+}(\mu),\ldots , J_{n}^{+}(\mu)\}\hbox{ and
}{\cal F}^{-}(\mu)=\{J_{1}^{-}(\mu),\ldots , J_{n}^{-}(\mu)\}.\hskip.25in
\end{equation}
The elements of ${\cal F}^{+}(\mu)$ will be called {\it right intervals} and the elements of ${\cal F}^{-}(\mu)$ will be called {\it
left intervals}.

\demo{{R}emark} Most of our results and definitions will be given for
right intervals only. The corresponding facts for left intervals can be easily
obtained in a symmetrical way or by applying the given ones to the reflected
measure $\tilde{\mu}=\sum_{i=1}^{n}k_{i}\delta_{-y_{i}}$.
\enddemo

The role of the gap interval in the covering properties of the $I_{i,j}$'s can
be seen by the following (see \cite{Mel}):

\specialnumber{2}\proclaim{Proposition}
{\rm (i)} Let $1\leq i<j\leq n$. Then $I_{i,j}\neq\emptyset$ if and only if
$J_{i}^{+}\cap J_{j}^{-}\neq\emptyset$.
\vglue4pt
{\rm (ii)} If $a_{j}\notin J_{i}^{+}$ and $a_{i}\notin J_{j}^{-}$ then $\left|
I_{i,j}\right|  =\left|  J_{i}^{+}\cap J_{j}^{-}\right|  $.
\vglue4pt
{\rm (iii)} If $\mu$ is admissible then $\left|  J(\mu)\right|  =\left|
G(\mu)\right|  $ and $J(\mu)\subseteq J_{1}\cup\cdots\cup  J_{n}$.
\endproclaim

Any interval $I_{i,j}$ as in Proposition 2(ii) will be called {\it special}. We also have the following.

\specialnumber{2}\proclaim{Lemma}
The interval $I_{i,j}\neq\emptyset$ is special if and only if $\left|
I_{i,j}\right|  <\min(k_{i},k_{j})$.
\endproclaim

\demo{Proof}
It is easy to see that $\left|  I_{i,j}\right|  =\max(k_{i}+k_{j}-(a_{j}%
-a_{i}),0)$. Hence if nonempty it would be special if and only if $a_{j}%
>a_{i}+k_{i}$ and $a_{i}<a_{j}-k_{j}$ and this easily completes the proof.
\enddemo

To proceed further for each fixed $i$ we set $l_{i}=\min\{l\leq i:a_{l}\in
J_{i}^{-}\}$, $r_{i}=\max\{r\geq i:a_{r}\in J_{i}^{+}\}$ and define the intervals%
\begin{equation}
F_{i}=F_{i}(\mu)=[y_{i}-K_{l_{i}}^{i},\,y_{i}+K_{i}^{r_{i}}].%
\end{equation}
Then the following holds (see \cite{Mel}). 

\specialnumber{3}\proclaim{Proposition}
{\rm (i)}We have $F_{i}=I_{l_{i},i}\cup I_{i,l_{i}+1}\cup\cdots\cup  I_{i,i}\cup
I_{i,i+1}\cup\cdots\cup  I_{i,r_{i}}$.
\vglue4pt 
{\rm (ii)} For any $i$ the nonempty of the closed intervals $I_{1,i},\ldots , I_{l_{i}%
-1,i}$ and $I_{i,r_{i}+1},\ldots , I_{i,n}$ {\rm (}\/if any\/{\rm )} are pairwise disjoint and each
of them is disjoint from $F_{i}$.

\vglue4pt 
{\rm (iii)} The set $E(\mu)$ is covered by the $n$ main intervals $F_{i}${\rm ,} $1\leq
i\leq n$ together with the nonempty {\rm (}\/if any\/{\rm )} special intervals $I_{p,q}$ where
$a_{q}\notin J_{p}^{+}$ and $a_{p}\notin J_{q}^{-}$.
\endproclaim

By exploiting the above structure of the gap interval we will prove the
following basic for our developments (see also \cite{Mel1}).

\specialnumber{4}\proclaim{Proposition}
{\rm (i)} The set $G(\mu)$ can be covered by appropriately placing certain parts of
the nonempty of the intervals $J_{i}^{+}\cap J_{j}^{-}$ over $[y_{i}%
+k_{i},y_{j}-k_{j}]$ for $1\leq i<j\leq n${\rm ,} each such part used at most once.

\vglue4pt 
{\rm (ii)} In particular if $\mu$ is admissible $J(\mu)$ can be also covered as in
{\rm (i),} where each used part of $J_{i}^{+}\cap J_{j}^{-}$ is placed appropriately
over $[a_{i},a_{j}]$.
\endproclaim

\demo{Proof}
(i) Consider an $i\,$\ with $1\leq i\leq n$. If $a_{i}\notin J_{s}$ for every
$l_{i}\leq s\leq r_{i}$ with $s\neq i$, then clearly $|  J_{i}^{+}\cap
J_{s}^{-}|  =k_{s}$ for any $i<s\leq r_{i}$ (respectively $|
J_{s}^{+}\cap J_{i}^{-}|  =k_{s}$ for any $l_{i}\leq s<i$) and so
writing $\tilde{I}_{i,s}=[y_{i}+K_{i}^{s-1},\,y_{i}+K_{i}^{s}]\break\subseteq
I_{i,s}$ (respectively $\tilde{I}_{s,i}=[y_{i}-K_{s}^{i},\,y_{i}-K_{s+1}%
^{i}]\subseteq I_{s,i}$) we \pagebreak easily conclude that these intervals cover
$F_{i}\backslash I_{i,i}$ and have lengths equal to $|  J_{i}^{+}\cap
J_{s}^{-}|  $ (respectively $|  J_{s}^{+}\cap J_{i}^{-}|  $)
and using (\ref{a3}) each such $\tilde{I}_{i,s}$ (respectively $\tilde
{I}_{s,i}$) is contained in $[y_{i},y_{s}]$ (respectively $[y_{s},y_{i}]$).
$\phantom{\sum^1}$

Now assume that there is a largest possible $s$ such that $i<s\leq r_{i}$ and
$a_{i}\in J_{s}^{-}$. Then since also $a_{s}\in J_{i}^{+}$ we conclude that
$[a_{i},a_{s}]=J_{i}^{+}\cap J_{s}^{-}$ and so the part of $G(\mu)$ that lies
in $[y_{i}+k_{i},y_{s}-k_{s}]$ can be obviously covered by using certain parts
of just $J_{i}^{+}\cap J_{s}^{-}$. The remaining part of the $F_{i}\cap
(y_{i},+\infty)$ that is $F_{i}\backslash(-\infty,y_{s}+k_{s})$ (if any) has
length $$(y_{i}+K_{i}^{r_{i}})-(y_{s}+k_{s})=K_{i}^{r_{i}}-(a_{s}-a_{i}%
+2K_{i}^{s} -k_{i})<K_{s+1}^{r_{i}}$$ and is thus covered by the intervals
$$\tilde{I}_{i,j}=[y_{i}+K_{i}^{j-1},\,y_{i}+K_{i}^{j}]\subseteq I_{i,j}$$ where
$s<j\leq r_{i}$ each contained in the corresponding $[y_{i},y_{j}]$ and having
length $\left|  J_{i}^{+}\cap J_{j}^{-}\right|  $ since $a_{i}\notin J_{j}$
for every $s<j\leq r_{i}$. Similar considerations can be applied if $a_{i}\in
J_{s}^{+}$ for some $l_{i}\leq s<i$.

Finally for any special interval $I_{p,q}$ where $a_{q}\notin J_{p}^{+}$ and
$a_{p}\notin J_{q}^{-}$ we know that $\left|  I_{p,q}\right|  =\left|
J_{p}^{+}\cap J_{q}^{-}\right|  $.

These, combined with Proposition 3(iii), complete the proof of (i), obsering
that any part of any used piece that is contained in $$\bigcup\limits_{i=1}%
^{n}I_{i,i}=\bigcup\limits_{i=1}^{n}[y_{i}-k_{i},y_{i}+k_{i}]$$ can be ignored.
\vglue4pt
(ii) If $\mu$ is admissible then all gaps in $[y_{1}-k_{1},y_{n}%
+k_{n}]\backslash(I_{1,1}\cup\cdots\cup  I_{n,n})$ are covered and so $\left|
G(\mu)\right|  =\left|  J(\mu)\right|  $. Therefore we can via the mapping
$Q^{-1}$ transport the way $G(\mu)$ is covered to cover $J(\mu)$ and this
completes the proof observing that any piece placed over $[y_{i}+k_{i}%
,y_{j}-k_{j}]$ when transported via $Q^{-1}$ will lie over $[a_{i},a_{j}]$.
\enddemo

{\it Remarks}.  (i) When the covering of $G(\mu)$ that is described in the
above proof is transported via $Q^{-1}$ to cover $J(\mu)$ some intervals might
shrink due to existence of intermediate masses. Here the fact that $Q^{-1}$ is
distance nonincreasing is used.
\vglue4pt
(ii) It is evident from the proof of Proposition 4 that in the case $a_{j}\in
J_{i}^{+}$ and $a_{i}\in J_{j}^{-}$ the whole part $[a_{i},a_{j}]$ of the gap
interval is equal and hence completely covered by $J_{i}^{+}\cap J_{j}^{-}$.
However \pagebreak due to the possible existence of masses between $y_{i}$ and
$y_{j}$, it might be necessary to break $J_{i}^{+}\cap J_{j}^{-}$ into several pieces
before placing it over $[y_{i}+k_{i},y_{j}-k_{j}]$. Actually this is the only
case where such a breaking occurs.
\vglue8pt

It would be important to keep track of exactly how the parts of the $J_{i}%
^{+}\cap J^{-}_{\raise4pt\hbox{$\scriptstyle j$}}$'s are placed to cover $G(\mu)$ and
$J(\mu)$. This has been more or less analysed in the 
above proof except for the case of special intervals. Related to this we have the
following (where by
$\frak{l}(I)$,
$\frak{r}(I)$ we will denote the left and right endpoints of the interval $I$).

\specialnumber{3}\proclaim{Lemma}
Suppose that $1\leq i\leq n${\rm ,} that $r_{i}\leq r<s$ and that both $I_{i,r}$ and
$I_{i,s}$ are nonempty. Then
\begin{equation}
\frak{l}(I_{i,s})-\frak{r}(I_{i,r})=\mathop{\rm dist}\nolimits (a_{s},J_{i}%
)+K_{r+1}^{s-1}%
\end{equation}
and a similar relation holds when $s<r\leq l_{i}$.
\endproclaim

\demo{Proof}
We have $\frak{l}(I_{i,s})-\frak{r}(I_{i,r})=(y_{s}-K_{i}^{s})-(y_{i}%
+K_{i}^{r})$ and using the relation $y_{s}-y_{i}=a_{s}-a_{i}+k_{i}%
+2k_{i+1}+\cdots + 2k_{s-1}+k_{s}$ we easily get $\frak{l}(I_{i,s})-\frak{r}%
(I_{i,r})=a_{s}-a_{i}-k_{i}+k_{r+1}+\cdots + k_{s-1}=a_{s}-\frak{r}(J_{i}%
)+K_{r+1}^{s-1}$ which completes the proof since $a_{s}>a_{i}$ and
$a_{s}\notin J_{i}$.
\enddemo

{\it Remarks}. (i) Clearly $\frak{l}(I_{i,r})=\frak{l}(F_{i})$ if $r=r_{i}
$. Thus Lemma 3 shows where the special intervals are located after the
related $F_{i}$'s. For example it shows that there is a gap between $F_{i}$
and the first special interval of the form $I_{i,s}$ (if any) that is at least
$\mathop{\rm dist}\nolimits (a_{s},J_{i})$ and in case $\mu$ is admissible has to be
covered by intervals of the form $I_{p,q}$ where $p\neq i$ and $q\neq i$. This
exact location will be important in our proof of Theorem 1.
\vglue4pt 
(ii) Actually the above results show how one can read off the covering
properties of the family of intervals $I_{i,j}(\mu)$ for $i<j$ from the
corresponding overlappings of the families ${\cal F}^{+}(\mu)$ and
${\cal F}^{-}(\mu)$ over the gap interval. In particular they show that the
length and exact location in $E(\mu)$ of the special intervals $I_{i,r}$ (if
any) depend only on the behavior of the gap interval and the corresponding
$J_{m}^{-}$'s that are located to the right of the {\it right endpoint }of
$J_{i}^{+}$.

\demo{Notation} (i) In this paper we will use the notation $\left|
\cdots\right|  $ in two different contexts: If $S$ is a subset of ${\Bbb  R}$
(which will ususaly be the union of finitely many closed intervals) then
$\left|  S\right|  $ will denote its Lebesgue measure. If on the other hand
$T$ is a finite set (that will usually consist of a finite number of intervals)
then $\left|  T\right|  $ will denote the cardinality of $T$.
\vglue4pt
(ii) For every family ${\cal U}$ of intervals by $\bigcup{\cal U}$ we
will denote the union of all elements of ${\cal U}$.
\vglue4pt
(iii) As above for any interval $I\subseteq{\Bbb  R}$ by $\frak{l}(I)$,
$\frak{r}(I)$ we will denote its left and right endpoints respectively.
\enddemo
 
\section{The measure $\mu$}

Let
\begin{equation}
\gamma=\dfrac{-1+\sqrt{61}}{12}=0.5675208\ldots
\end{equation}
be the positive solution of the quadratic equation
\begin{equation}
12\gamma^{2}+2\gamma-5=0.\label{eq}%
\end{equation}

Assuming that $C>1+\gamma$ there must exist measures $\mu$ as in (\ref{a1})
such that $R(\mu)>1+\gamma$. We then consider the smallest possible integer
$n$ such that there exists a measure $\mu=\sum_{i=1}^{n}k_{i}\delta_{y_{i}}$
such that $R(\mu)>1+\gamma$. Then $R(\nu)\leq1+\gamma$ for any measure as in
(\ref{a1}) that contains less than $n$ positions. Hence using Lemma 1 and
Proposition 1(ii) we may assume that $\mu$ is admissible; that is, it satisfies
(\ref{a2}) and (\ref{a3}).

Moreover we may assume that all the $y_{i}$'s and all the $k_{i}$'s are
positive {\it integers}. Indeed we can find rational numbers $k_{i}%
^{\prime}>k_{i}$ and $y_{i}^{\prime}$ for $1\leq i\leq n$ such that
$0<y_{i+1}^{\prime}-y_{i}^{\prime}<y_{i+1}-y_{i}$, the $y_{i}^{\prime}$ and
$k_{i}^{\prime}$ satisfy (\ref{a2}) and the (as it is easy to see) admissible
measure $\mu^{\prime}=\sum_{i=1}^{n}k_{i}^{\prime}\delta_{y_{i}}^{\prime}$
still satisfies $R(\mu^{\prime})>1+\gamma$. Then by multiplying all
$y_{i}^{\prime}$ and $k_{i}^{\prime}$ by an appropriate integer we get a
measure with all entries integers.

From now on we will fix such a measure $\mu$ and let its gap interval $J(\mu)$
and its corresponding cover ${\cal F}(\mu)=\{J_{1},\ldots , J_{n}\}$ be as in
Section 2.

Then we write
\begin{equation}
J(\mu)=[0,N]=\omega_{1}\cup\cdots\cup \omega_{N},
\end{equation}
where $N$ is a positive integer and $\omega_{p}=[p-1,p]$ for $p=1,2,\ldots , N$.
Each $\omega_{p}$ will be called a {\it place} in the gap interval $J(\mu
)$. Also since the corresponding $x_{i}$ and $k_{i}$'s are integers to each
such $\omega_{p}$ there correspond three nonnegative integers $h_{p}^{+}$,
$h_{p}^{-}$ and $h_{p}$ such that
\begin{equation}
h_{p}^{+}=\sum_{i=1}^{n}\chi_{J_{i}^{+}}(x)\hbox{, }h_{p}^{-}=\sum_{i=1}%
^{n}\chi_{J_{i}^{-}}(x)\hbox{ and }h_{p}=h_{p}^{+}+h_{p}^{-}%
\end{equation}
for any $x\in {\rm int} (\omega_{p})$. Clearly
\begin{equation}
2K_{1}^{n}=\sum_{i=1}^{n}\left|  J_{i}\right|  \geq h_{1}+\cdots + h_{N}.\label{br1}%
\end{equation}
(We write $\geq$ since $J_{1}\cup\cdots\cup  J_{n}$ might contain points outside
$J(\mu)$.)

We will be considering that over each place $\omega_{p}$ there are $h_{p}$
distinct intervals of length $1$ which we call {\it bricks }$h_{p}^{+}$
corresponding to the right intervals that contain $\omega_{p}$ and $h_{p}^{-}$
to the left. It is clear that $h_{1}+\cdots + h_{N}$ is the total number of bricks.

We also let
\begin{equation}
P=\{a_{1},\ldots , a_{n}\}
\end{equation}
denote the set of all {\it positions} (centers of the $J_{i}$'s) in the gap interval.

Now we consider the set of places
\begin{equation}
E_{1}=\{\omega_{p}\subseteq J(\mu):h_{p}=1\}
\end{equation}
over which exactly one interval from the family ${\cal F}^{+}(\mu
)\cup{\cal F}^{-}(\mu)$ passes. It is then easy to see, using (\ref{br1})
and Proposition 2(iii) that the places in $E_{1}$ are the only ones that have
the property of pushing $R(\mu)$ to something bigger than $\frac{1}{2}$. Thus
it would be important to analyze the behavior of the intervals of
${\cal F}^{+}(\mu)\cup{\cal F}^{-}(\mu)$ that contain such places. We
will consider only right intervals the corresponding statements for left
intervals being symmetrical. It is clear, by Proposition 4(ii), that if a
$J_{i}^{+}$ contains an $\omega_{p}\in E_{1}$ then $\omega_{p}$ can be covered
only through the involvement of \ this $J_{i}^{+}$.

There are essentially two cases to consider. The first is treated in the following.

\specialnumber{5}\proclaim{Proposition}
Suppose that for some $i\geq1$ there exist $\omega_{p}\in E_{1}$ and $x\in
 {\rm int} (\omega_{p})\subseteq J_{i}^{+}$ such that $Q(x)\leq\frak{r}(F_{i})$.
Then we have
\begin{equation}
(a_{i},x]\cap P=\emptyset\label{a10}%
\end{equation}
and
\begin{equation}
a_{i+1}-a_{i}\leq K_{i+1}^{r_{i}}=\left|  F_{i}\backslash(-\infty,y_{i}%
+k_{i})\right|  .\label{a11}%
\end{equation}
\endproclaim

\demo{Proof}
Suppose that $(a_{i},x]\cap P=\{a_{i+1},\ldots , a_{s}\}\neq\emptyset$ and so
$a_{s}\leq x<a_{s+1}$. Since $h_{p}=1$ it is clear that no interval other than
$J_{i}^{+}$ contains $x$ and so by Proposition 2(i) we have $I_{s,r}%
=\emptyset$ whenever $r>s$. Hence moving $k_{s}\delta_{y_{s}}$ to the left
by $a_{s}-a_{s-1}$ will not change the connectivity of $E(\mu)$ since this
mass does not  interact with any mass to its left, since the inequality
$Q(x)\leq\frak{r}(F_{i})$ implies that $y_{s}$ belongs to $F_{i}$ that will
hence not change, as long as $a_{s}\in J_{i}^{+}$, and since this movement can
only enlarge the intervals $I_{l,s}$ for $l<i$. But then the resulting measure
$\mu^{\prime}$ will have the same $E(\mu^{\prime})$ but will not satisfy the
separability condition (\ref{a2}) for the $s-1$ position. However in view of
Lemma 1 this implies that there is a measure $\mu^{\prime\prime}$ containing
at most $n-1$ positions with $R(\mu^{\prime\prime})\geq R(\mu^{\prime}%
)=R(\mu)$ and this contradicts our choice of $\mu$. Hence $(a_{i},x]\cap
P=\emptyset$.

Next we will show that $\frak{r}(F_{i})<y_{i+1}-k_{i+1}$ is impossible. Indeed
if this happened then since $x<a_{i+1}$ it is easy to see that $I_{l,s}%
=\emptyset$ whenever $l<i<s$ and so the  \pagebreak interval $[\frak{r}(F_{i}%
),\frak{r}(F_{i})+1]$ must be covered by some $I_{i,s}$ where necessarily
$s>r_{i}$ and so $I_{i,s}$ is a special interval. Thus $\frak{l}(I_{i,s}%
)\leq\frak{r}(F_{i})$ which contradicts Lemma 3. Hence $\frak{r}(F_{i}%
)=y_{i}+K_{i}^{r_{i}}\geq y_{i+1}-k_{i+1}$ and since $y_{i+1}-y_{i}%
=a_{i+1}-a_{i}+k_{i+1}+k_{i}$ we get (\ref{a11}).
\enddemo

If for the right interval $J_{i}^{+}$ there exist $\omega_{p}\in E_{1}$ and
$x$ $\in {\rm int} (\omega_{p})\subseteq J_{i}^{+}$ such that $Q(x)\leq
\frak{r}(F_{i})$ (and so $(a_{i},x]\cap P=\emptyset$) then the right
interval $J_{i}^{+}$ will be called {\it clean}. A symmetrical definition
applies for the left intervals $J_{j}^{-}$.

Suppose now that for some $m\geq1$ the right interval $J_{m}^{+}$ contains at
least one place from $E_{1}$ but is not clean. Then defining
\begin{equation}
w=\min\{q:\omega_{q+1}\subseteq J_{m}^{+}\hbox{ and\ }h_{q+1}=1\}\geq a_{m}%
\end{equation}
we must have $(a_{m},w]\cap P\neq\emptyset$. Indeed if $(a_{m},w]\cap
P=\emptyset$ then clearly $w+1\leq a_{m+1}$ and moreover since $[w,w+1]\in
E_{1}$ the {\it interval} $Q((a_{m},w+1])\subseteq\lbrack y_{m},y_{m+1}]$
must be covered only by intervals of the form $I_{m,r}$ for $r>m$ (because by
Proposition 2(i), $I_{l,r}=\emptyset$ whenever $l<m<m+1\leq r$). However
Proposition 3(ii) now implies that we must have $Q((a_{m},w+1])\subseteq
F_{m}$ and so $Q(w+\frac{1}{2})<\frak{r}(F_{m})$, which contradicts the assumption
that $J_{m}^{+}$ is not clean. Hence we may write
\begin{equation}
(a_{m},w]\cap P=\{a_{m+1},\ldots , a_{s}\}\neq\emptyset.
\end{equation}
Clearly $h_{p}\geq2$ for all $a_{m}\leq p\leq w$. Now let
\begin{equation}
g(J_{m}^{+})=a_{s}-a_{m}\hbox{, }K(J_{m}^{+})=K_{m+1}^{s}.%
\end{equation}
Then we have the following.

\specialnumber{4}\proclaim{Lemma}
The interval $(y_{s}+k_{s},y_{s}+k_{s}+1]$ must be covered by a special
interval $I_{m,t}$ for some $t>r_{m}$. Moreover we must have
\begin{equation}
g(J_{m}^{+})+K(J_{m}^{+})\geq\mathop{\rm dist}\nolimits (a_{t},J_{m})+K_{s+1}%
^{t-1}.\label{m3}%
\end{equation}
\endproclaim

\demo{Proof}
By a similar reasoning as in the proof of Proposition 5, \ we conclude that
$F_{m}$ cannot cover the point $y_{s}+k_{s}+\frac{1}{2}$. Since for any
$l\leq s<r$ we have $I_{l,r}=\emptyset$ unless $l=m$ we conclude that it
must be covered by some special interval $I_{m,t}$ for some $t>r_{m}$ and so
$a_{t}>a_{m}+k_{m}=\frak{r}(J_{m})$. Since the $y_{l}$'s and the $k_{l}$'s are
integers we have
\begin{equation}
y_{s}+k_{s}\geq\frak{l}(I_{m,t})=y_{t}-k_{m}-\cdots-k_{t}.\label{m4}%
\end{equation}
Writing now $$y_{s}+k_{s}=y_{m}+a_{s}-a_{m}+k_{m}+2k_{m+1}+\cdots + 2k_{s}$$
 and
$$y_{t}-k_{m}-\cdots-k_{t}=y_{m}+a_{t}-a_{m}+k_{m+1}+\cdots + k_{t-1}$$ we get
(\ref{m3}).
\enddemo

{\it Remark}. In the above lemma we may actually assume that equality holds
in (\ref{m4}) and hence also in (\ref{m3}). Indeed clearly the mass
$k_{s}\delta_{y_{s}}$ interacts with no mass to the right of it (meaning that
$I_{s,j}=\emptyset$ for every $j>s$). Hence as in the proof of Proposition 5
it can be moved to the left until either equality in (\ref{m4}) occurs or the
separability inequality (\ref{a2}) for $i=s-1$ is violated. But as in the
proof of that proposition the second alternative cannot happen.

\section{Further covering properties of $\mu$}

By Proposition 4 and since $\mu$ is admissible to each $\omega_{p}$ we can
associate an $\omega_{c(p)}$ and certain $i(p)<j(p)$ such that $\omega_{p}%
\in\lbrack a_{i(p)},a_{j(p)}]$, $\omega_{c(p)}\subseteq J_{i(p)}^{+}\cap
J_{j(p)}^{-}$ and such that the part $\omega_{c(p)}$ of $J_{i(p)}^{+}\cap
J_{j(p)}^{-}$ is used (corresponds to the part of $I_{i(p),j(p)}$ used) to
cover $\omega_{p}\subseteq J(\mu)$ (equivalently $Q(\omega_{p})\subseteq
\lbrack y_{1},y_{n}]$) according to above mentioned proposition. Moreover it
is clear that the mapping
\begin{equation}
p\rightarrow(c(p),i(p),j(p))\label{m1}%
\end{equation}
is {\it one\/{\rm -}\/to\/{\rm -}\/one}. We will write $\omega_{c(p)}\rightarrow\omega_{p}$ and
we will say that that $\omega_{c(p)}$ {\it covers }$\omega_{p}$. Also to
indicate the exact way this covering takes place we will say that $\omega_{p}$
is covered by $(\omega_{c(p)},J_{i(p)}^{+},J_{j(p)}^{-})$ and we will say that
$\omega_{p}$ is covered by $\omega_{c(p)}$ through the {\it interaction} of
the right interval $J_{i(p)}^{+}$ with the left interval $J_{j(p)}^{-}$.

\demo{{R}emark} It may happen that $\omega_{p}$ is covered by more than one
way according to Proposition 4. In such a case we choose exactly one of these
ways arbitrarily to make the mapping $c$ well defined.
\enddemo

For any $\omega_{p}$ that covers at least one place we let
\begin{equation}
l(p)=\min\{i:\omega_{p}\subseteq J_{i}^{+}\}<r(p)=\max\{j:\omega_{p}\subseteq
J_{j}^{-}\}
\end{equation}
(both well defined) and we define the intervals
\begin{equation}
L_{p}=J_{l(p)}^{+}\hbox{ and }R_{p}=J_{r(p)}^{-}.%
\end{equation}

Now except for $E_{1}$ we will more generally consider for any nonnegative
integers $s,t$ the sets
\begin{equation}
E_{s,t}=\{\omega_{p}\subseteq J(\mu):h_{p}^{+}=s\hbox{ and }h_{p}^{-}=t\}
\end{equation}
and
\begin{equation}
E_{t}=\{\omega_{p}\subseteq J(\mu):h_{p}=t\}=\bigcup_{a+b=t}E_{a,b}.%
\end{equation}

We have the following.

\specialnumber{5}\proclaim{Lemma}
{\rm (i)} $\omega_{p}\in E_{a,b}$ can cover at most $a.b$ places in $J(\mu)$.
\vglue4pt
{\rm (ii)} Any $\omega_{p}$ can cover at most $h_{p}-1$ places in $E_{1}\cup
E_{1,1}$.
\endproclaim

\demo{Proof}
For (i) obviously $a.b$ is equal to the number of all possible pairs $(A,B)$
of a right interval $A$ and a left interval $B$ such that $\omega_{p}\subseteq
A\cap B$. We will now prove (ii). If $\omega_{p}$ covers at least one place
then $l(p),r(p)$ are well defined. Suppose that for some $i,j $ with
$l(p)<i<j<r(p)$ a place $\omega_{q}\in E_{1}\cup E_{1,1}$ is covered through
$(\omega_{p},J_{i}^{+},J_{j}^{-})$. Then we have $\omega_{q}\subseteq\lbrack
a_{i},a_{j}]$. However $\omega_{p}\subseteq J_{i}^{+}\cap J_{j}^{-}$ so it is
clear that $\chi_{J_{i}^{+}}+\chi_{J_{l(p)}^{+}}\geq2$ on $[a_{i},p]$ and
$\chi_{J_{j}^{-}}+\chi_{J_{r(p)}^{-}}\geq2$ on $[p-1,a_{j}]$. Therefore
$h_{q}^{+}\geq2$ if $q\leq p$ and $h_{q}^{-}\geq2$ if $q\geq p-1$ and both lead 
to a contradiction. Hence the possible $\omega_{q}^{\phantom{1}}\in E_{1}\cup
E_{1,1}$ covered by $\omega_{p}$ can come only from interactions in which at least
one of the intervals $L_{p}$ and $R_{p}$ is involved and it easy to see that there
are $(h_{p}^{+}-1)+(h_{p}^{-}-1)+1=h_{p}-1$ such interactions.
\enddemo

{\it Remark}. This lemma in particular implies that an $\omega_{p}$ in
$E_{1}$ does not cover any place, an $\omega_{p}$ in $E_{2}$ covers at most one
place (and this can happen only if $h_{p}^{+}=h_{p}^{-}=1$) and an $\omega
_{p}$ in $E_{3}$ covers at most two places. Also an $\omega_{p}\in E_{3,1}\cup
E_{1,3}$ can cover at most three places whereas an $\omega_{p}\in E_{2,2}$ can
cover at most four places at most three of which can belong to $E_{1}\cup
E_{1,1}$.
\vglue4pt

We will introduce now the following notation: Suppose, for example, that an
$\omega_{p}\in E_{3}$ covers an $\omega_{q}\in E_{1}$ and also an $\omega
_{a}\in E_{1,1}\subseteq E_{2}$ that in turn covers an $\omega_{b}\in E_{1}$.
Then we will say that $\omega_{p}$ is the {\it head} of an $E_{3}%
\rightarrow(E_{1},(E_{2}\rightarrow E_{1}))$ {\it pattern}. We will
consider the following nine types of such patterns:
$$
\begin{array}{lcl}
\hbox{Type }1\hbox{{}}  & :&\hbox{ }E_{1}\\
\hbox{Type }2\hbox{{}}  & :&\hbox{ }E_{2}\rightarrow E_{1}\\
\hbox{Type }3\hbox{{}}  & :&\hbox{ }E_{2}\rightarrow E_{2}\rightarrow E_{1}\\
\hbox{Type }4\hbox{{}}  & :&\hbox{ }E_{2}\rightarrow E_{2}\rightarrow
E_{2}\rightarrow E_{1}\\
\hbox{Type }5\hbox{{}}  & :&\hbox{ }E_{2}\rightarrow(E_{3}\rightarrow
(E_{1},E_{1}))\\
\hbox{Type }6\hbox{{}}  & :&\hbox{ }E_{3}\rightarrow(E_{1},E_{1})\\
\hbox{Type }7\hbox{{}}  & :&\hbox{ }E_{3}\rightarrow((E_{2}\rightarrow
E_{1}),E_{1})\\
\hbox{Type }8\hbox{{}}  & :&\hbox{ }E_{1,3}\cup E_{3,1}\rightarrow(E_{1}%
,E_{1},E_{1})\\
\hbox{Type }9\hbox{{}}  & :&\hbox{ }E_{4}\rightarrow((E_{3}\rightarrow
(E_{1},E_{1})),(E_{2}\rightarrow E_{1}),E_{1},E_{1}).%
\end{array}
$$
It is required that the $E_{1}$'s appearing in the Types $5,6,8$ and $9$
patterns are referring to distinct places. It is also clear that if $\omega
_{p}$ is the head of a Type $j$ pattern then for $1\leq j\leq5$ we must have
$\omega_{p}\in E_{1,1}$ and for $j=6,7$ we must have $\omega_{p}\in
E_{1,2}\cup E_{2,1}$. The possibility $\omega_{p}\in E_{2,2}$ has been
excluded from the Type $8$ pattern.  

Moreover we have the following.

\specialnumber{6}\proclaim{Lemma}
 Consider any Type $j$ pattern where $1\leq j\leq9$ and let $T$ be the
set of all places involved in it. Then\/{\rm :}\/
\begin{itemize}
\ritem{(i)} All places indicated in this pattern are distinct{\rm ;} hence $T$ has as many
elements as the $E_{t}$\/{\rm '}\/s appearing in the pattern.

\ritem{(ii)} No $\omega_{q}\in T$ can cover any place outside $T$.

\ritem{(iii)} If an $\omega_{q}$ covers the head of this pattern{\rm ,} then $\omega
_{q}\notin T$.

\ritem{(iv)} Given $\omega_{q}\in T$ and a pair $(A,B)$ of a right interval $A$ and a
left interval $B$ such that $\omega_{q}\subseteq A\cap B$ then there exists
$\omega_{s}\in T$ such that $(\omega_{q},A,B)$ covers $\omega_{s}$.
\end{itemize}

\endproclaim

\demo{Proof}
For (i) it obviously suffices to consider only places in the same $E_{t}$ that
are covered by places in the same $E_{s}$. Hence by the requirements set for
the Types $5,6,8$ and $9$ it only remains to treat the Types $3$ and $4$.
Suppose for example that a Type $4$ pattern involves $\omega_{a}%
\rightarrow\omega_{b}\rightarrow\omega_{p}\rightarrow\omega_{q}$ but
$\omega_{a}=\omega_{p}$. Then $\omega_{a}\in E_{2}$ would have to cover the
two different places $\omega_{b}\in E_{2}$ and $\omega_{q}\in E_{1}$
contradicting Lemma 5 The proof for the other cases is similar. The assertion
(ii) follows again by Lemma 5, (iii) can be proved in a similar way as (i) and
(iv) can be proved by examining each considered pattern.
\enddemo

Let $u_{j}$ denote the number of places in a Type $j$ pattern and $v_{j}$ the
corresponding number of bricks. Then clearly $u_{1}=v_{1}=1$, $u_{2}=2$,
$v_{2}=3$, $u_{3}=u_{6}=3$, $v_{3}=v_{6}=5$, $u_{4}=u_{5}=u_{7}=u_{8}=4$,
$v_{4}=v_{5}=v_{7}=v_{8}=7$, $u_{9}=8$ and $v_{9}=14$. Also for $1\leq j\leq9
$ let
\begin{equation}
\lambda_{j}=u_{j}-\gamma v_{j}.%
\end{equation}
It is easy to see that
\begin{equation}
0<\lambda_{4}=\lambda_{5}=\lambda_{7}=\lambda_{8}<\lambda_{9}<\lambda
_{3}=\lambda_{6}<\lambda_{2}<\lambda_{1}.\label{l1}%
\end{equation}

Now for any $\omega_{p}$ that is not the head of any Type $j$ pattern for any
$1\leq j\leq9$ we let $T_{p}$ be the set that consists of $\omega_{p}$ and all
places from all (maximal) patterns whose head is covered by $\omega_{p}$ and
let
\begin{equation}
H_{p}=\sum_{\omega_{s}\in T_{p}}h_{s}%
\end{equation}
be the corresponding number of bricks that lie over all such places.

If now $\omega_{p}$ is the head of a Type $j$ pattern for some $1\leq j\leq9$
we let $T_{p}$ be the set of all places involved in this pattern, so $\left|
T_{p}\right|  =u_{j}$, but let
\begin{equation}
H_{p}=v_{j}+1\label{h1a}%
\end{equation}
in this case (instead of $v_{j}$). This modification, whose use will be made
clear later, results in the following estimate
\begin{equation}
\left|  T_{p}\right|  \leq\dfrac{8}{15}H_{p}<\gamma H_{p}%
\end{equation}
whenever $\omega_{p}$ is the head of such a pattern.

We also define $T_{p}=\emptyset$ and $H_{p}=0$ if $\omega_{p}$ does not fall
into one of the above two categories (for example an $\omega_{p}\in E_{2}$
that say covers an $\omega_{q}\in E_{4}$).

We now have the following.

\specialnumber{7}\proclaim{Lemma}
For any $p\neq q$ the sets $T_{p}$ and $T_{q}$ {\rm (}\/if defined\/{\rm )} are either
disjoint or one of them is contained in the other.
\endproclaim

\demo{Proof}
We will associate to each $\omega_{s}\in T_{p}$ an integer $r=r(s)$, called
its rank, to be the length of the chain $\omega_{p}\rightarrow\cdots\rightarrow
\omega_{s}$ that leads to $\omega_{s}$. This is well defined since Lemma 6
implies that exactly one such chain can exist. Then if $T_{p}\cap T_{q}$ were
nonempty we choose an $\omega_{s}\in T_{p}\cap T_{q}$ whose rank in $T_{p}$ is
as small as possible. It is then clear that $\omega_{c(s)}$ cannot be
contained in both $T_{p}$ and $T_{q}$. Suppose that $\omega_{c(s)}\notin
T_{p}$ (the argument will show that the other case is impossible by the choice
of $\omega_{s}$). Then $\omega_{s}$ cannot be contained in any Type $j$
pattern whose head is covered by $\omega_{q}$ since this would easily imply
that $\omega_{c(s)}$ is either contained in the same pattern or is equal to
$\omega_{q}$ and in both cases $\omega_{c(s)}\in T_{p}$. The only alternative
is that $\omega_{s}=\omega_{q}$ and so that $\omega_{q}\in T_{p}$ must be the
head of a Type $j$ pattern. This easily implies that $T_{q}\subseteq T_{p}$
and completes the proof.
\enddemo

In the next two propositions we will show that any set $T_{p}$ will not
contribute significally to $R(\mu)>2(1+\gamma)$ unless $L_{p}$ and $R_{p}$
satisfy certain strong restrictions in relation with the set $E_{1}$.

\specialnumber{6}\proclaim{Proposition}
If $\omega_{p}$ is not the head of a Type $j$ pattern for any $1\leq j\leq9 $
and is such that at least one of the intervals $L_{p}$ and $R_{p}$ does not
contain any place from $E_{1}${\rm ,} then we have
\begin{equation}
\left|  T_{p}\right|  <\gamma H_{p}.\label{e1}%
\end{equation}
\endproclaim

\demo{Proof}
We may assume that $R_{p}$ does not contain any place from $E_{1}$, the proof
for $L_{p}$ being symmetrical. Let $h_{p}^{+}=a+1$ and $h_{p}^{-}=b+1$ and
number the the right intervals containing $\omega_{p}$ as $A_{0}=L_{p}%
,A_{1},\ldots , A_{a}$ and the left intervals containing $\omega_{p}$ as
$B_{0}=R_{p},B_{1},\ldots , B_{b}$ so that
\begin{equation}
\frak{l}(A_{0})<\frak{l}(A_{1})<\cdots <\frak{l}(A_{a})\hbox{ and }\frak{r}%
(B_{0})>\frak{r}(B_{1})>\cdots>\frak{r}(B_{b}).\hskip.4in
\end{equation}

Suppose first that $a,b>0$. By Lemma 5(ii), $\omega_{p}$ can cover the head of
a Type $j$ pattern with $1\leq j\leq5$ only if $A_{0}$ or $B_{0}$ is involved
(of course other patterns could also be so covered). However since
$\chi_{A_{0}}+\chi_{A_{1}}+\chi_{B_{0}}\geq2$ on $[\frak{l}(A_{1}%
),\min(\frak{r}(A_{1}),\frak{r}(B_{0}))]$ the triples $(\omega_{p},A_{i}%
,B_{0})$ for $i\geq1$ cannot cover an $E_{1}$ (since it should be contained in
$B_{0} $). Also since for any $\omega_{q}$ that is the head of a Type $6,7$ or
$9$ pattern there are exactly {\it two} intervals of the same direction
that contain it we conclude, using a similar argument as in the proof of Lemma
5, that $\omega_{p}$ can cover the head of such a pattern only if at least one
of the intervals $A_{0},A_{1},B_{0},B_{1}$ is involved. However if $i\geq2$
(so $b>1$) and $(\omega_{p},A_{1},B_{i})$ covers the head $\omega_{q}$ of a
Type $6$ pattern then we must have $\omega_{q}\in A_{1}\backslash B_{2}$ (and
so $q<p$) since $h_{q}^{+}\geq2$ if $\frak{l}(A_{1})\leq q\leq p$ and
$h_{q}^{-}\geq3$ if $p-1\leq q\leq\frak{r}(B_{2})$. Therefore $\omega_{q}$
would be contained in $A_{0}$ and $A_{1}$ and in exactly one other interval
$J$ of the opposite direction and moreover $(\omega_{p},A_{1},J)$ must cover a
place in $E_{1}$. But since $B_{0}$ doesn't contain places from $E_{1}$ we
clearly must have $\frak{r}(J)>\frak{r}(B_{0})$ and since $q<p$ this implies
that also $\omega_{p}\in J$. This contradicts the choice of $B_{0}=R_{p}$.
Hence $(\omega_{p},A_{1},B_{i})$ can cover only in Types $7,8$ or $9$.

Now similarly $\omega_{p}$ covers the head of a Type $8$ pattern only if at
least one of the intervals $A_{0},A_{1},A_{2},B_{0},B_{1},B_{2}$ is involved.
However if $i\geq2$ then $(\omega_{p},A_{i},B_{2})$ cannot cover the head of a
Type $9$ pattern since $h_{q}^{+}\geq3$ if $\frak{l}(A_{2})\leq q\leq p$ and
$h_{q}^{-}\geq3$ if $p-1\leq q\leq\frak{r}(B_{2})$. Also if $i\geq3$ then
$(\omega_{p},A_{2},B_{i})$ cannot cover the head of a Type $8$ (or $9$)
pattern for as before this would imply that this place must be in
$A_{2}\backslash B_{i}$ and this leads in a similar manner to a contradiction.

Hence the patterns covered by $\omega_{p}$ fall into exactly one of the
following categories:
\begin{itemize}
\item[(1)] With $A_{0}$ involved $\omega_{p}$ covers at most $b+1$ patterns of Type
 1--9.

\item[(2)] With $B_{0}$, but not $A_{0}$, involved $\omega_{p}$ covers at most $a$
patterns of Type 2--9.

\item[(3)] With $B_{1}$, but not $A_{0}$, involved $\omega_{p}$ covers at most $a$
patterns of Type 6--9.

\item[(4)] With $A_{1}$, but not $B_{0},B_{1}$, involved $\omega_{p}$ covers at most
$b-1$ patterns of Type 7--9.

\item[(5)] With $B_{2}$, but not $A_{0},A_{1}$, involved $\omega_{p}$ covers at most
$a-1$ patterns of Type $8$.
\end{itemize}

Let now $d_{i,j}$ the number of heads of Type $j$ patterns covered by
$\omega_{p}$ in the way described in category $(i)$ where $1\leq i\leq5,1\leq
j\leq9$. Some of those are of course $0$ as explained above, for example
$d_{4,6}=d_{5,9}=0$. Also we have given bounds for all five sums $\sum
_{j}d_{i,j}$, for example $\sum_{j}d_{4,j}\leq b-1$. Now it is clear that
\begin{equation}
\left|  T_{p}\right|  =1+\sum_{i,j}u_{j}d_{i,j}\hbox{ and }H_{p}%
=a+b+2+\sum_{i,j}v_{j}d_{i,j}.%
\end{equation}
Hence using (\ref{l1}) the bounds for the sums $\sum_{j}d_{i,j}$ and the zero
$d_{i,j}$'s we have
\begin{eqnarray*}
&&\hskip-18pt \left|  T_{p}\right|  -\gamma H_{p}=1+\sum_{i,j}\lambda_{j}d_{i,j}%
-\gamma(a+b+2) \\
&&\leq\ 1+\lambda_{1}\sum_{j}d_{1,j}+\lambda_{2}\sum_{j}d_{2,j}+\lambda_{6}%
\sum_{j}d_{3,j}+\lambda_{9}\sum_{j}d_{4,j} \\
&&\quad \ +\lambda_{7}\sum_{j}d_{5,j}-\gamma(a+b+2)\leq\hbox{ }(9-16\gamma
)(a+b)-(10-18\gamma)
\end{eqnarray*}
and so if $a+b\geq3$ we have
\begin{equation}
\left|  T_{p}\right|  -\gamma H_{p}\leq17-30\gamma<0.%
\end{equation}
If on the other hand $a=b=1$ and so $\omega_{p}\in E_{2,2}$ examining the five
categories it is easy to see that $\left|  T_{p}\right|  -\gamma H_{p}<0$
unless $d_{1,1}=2,d_{2,2}=d_{3,6}=1$ which implies that $\omega_{p}$ is the
head of a Type $9$ pattern, thus contradicting our assumption.

Suppose now that $a=0$ (the case $b=0$ is similar). Then $\omega_{p}$ covers
at most $b+1$ places and if $d_{j}$ of them are heads of Type $j$ patterns
then $\sum_{j}d_{j}\leq b+1$ and in a similar way we have
\begin{eqnarray*}
 \left|  T_{p}\right|  -\gamma
H_{p}&\hskip-8pt =\hskip-8pt&1+\sum_{j}\lambda_{j}d_{j}-\gamma(b+2)
\\ &\hskip-8pt=\hskip-8pt&\
1-\gamma-(2\gamma-1)(d_{1}+2d_{2}+3(d_{3}+d_{6})+4(d_{4}+d_{5}+d_{7}%
 +d_{8}))
\\ &\hskip-8pt&  -\ (15\gamma-8)d_{9}-\gamma\left(b+1-\sum_{j}d_{j}\right)
\end{eqnarray*}
and this would be negative unless $\sum_{j}d_{j}=b+1$ and
\[
d_{1}+2d_{2}+3(d_{3}+d_{6})+4(d_{4}+d_{5}+d_{7}+d_{8})+3.5d_{9}\leq3
\]
(and so $b\leq2$) since $\dfrac{15\gamma-8}{2\gamma-1}>3.5$, $\dfrac{1-\gamma
}{2\gamma-1}<3.3$ and the $d_{j}$'s are integers. These however easily imply
that $\omega_{p}$ must be the head of one of the Types 1--8 pattern which is
a contradiction. This completes the proof.
\enddemo

\specialnumber{7}\proclaim{Proposition}
If $\omega_{p}$ is not the head of a Type $j$ pattern for any $1\leq j\leq9 $
and is such that there is no $\omega_{s}\in L_{p}\cap R_{p}$ such that
$(\omega_{s},L_{p},R_{p})$ covers a place in $E_{1}${\rm ,} then we have
\begin{equation}
\left|  T_{p}\right|  <\gamma H_{p}.%
\end{equation}
\endproclaim

\demo{Proof}
By Propostion 6 both $L_{p}$ and $R_{p}$ contain places from $E_{1}$. Also by
the proof of that proposition we may assume that $h_{p}^{+}=a+1\geq2$ and
$h_{p}^{-}=b+1\geq2$. We number the the right and left intervals containing
$\omega_{p}$ as $A_{0}=L_{p},A_{1},\ldots , A_{a}$ and $B_{0}=R_{p},B_{1}%
,\ldots , B_{b}$ as in the proof of that proposition. By our assumption
$(\omega_{p},A_{0},B_{0})$ cannot cover the head of a Type $1$ pattern$.$

Suppose now that for some $i\geq1$, $(\omega_{p},A_{1},B_{i})$ covers the head
$\omega_{q}$ of a Type $j$ pattern for some $1\leq j\leq9$. If $\omega
_{q}\subseteq A_{1}\backslash B_{0}$ then clearly $h_{q}^{+}\geq2$ and so
$h_{q}\geq3$ and also there is no left interval $F$ such that $(\omega
_{q},A_{1},F)$ covers a place in $E_{1}$\ (since the only possible such $F$
would be $B_{0}$ which does not contain $\omega_{q}$). A similar statement
holds if $\omega_{q}\subseteq B_{1}\backslash A_{0}$. If $\omega_{q}\subseteq
A_{0}\cap B_{0}$ then also $h_{q}\geq3$ (since $\omega_{q}\subseteq A_{1}\cup
B_{1}$) and by our assumption $(\omega_{q},A_{0},B_{0})$ cannot cover any
place in $E_{1}$. Therefore the only possible values for $j$ are $7,8$ or $9$
and a similar statement holds if $(\omega_{p},A_{i},B_{1})$ covers the head
$\omega_{q}$ of a Type $j$ pattern.

Suppose now that for some $i\geq2$, $(\omega_{p},A_{2},B_{i})$ or $(\omega
_{p},A_{i},B_{2})$ covers the head $\omega_{q}$ of a Type $j$ pattern for some
$1\leq j\leq9$. Then $h_{q}^{+}\geq3$ or $h_{q}^{-}\geq3$ and so $j=8 $. If
$\omega_{q}\subseteq(A_{2}\backslash B_{0})\cup(B_{2}\backslash A_{0})$ then
as before it cannot happen that all places covered by $\omega_{q}$ are in
$E_{1}$, contradiction. Also if $\omega_{q}\subseteq A_{0}\cap B_{0}$ then
$(\omega_{q},A_{0},B_{0})$ cannot cover any place in $E_{1}$. Hence no such
covering can occur.

Therefore the patterns covered by $\omega_{p}$ fall into exactly one of the
following categories:
\begin{itemize}

\item[(1)] With $A_{0}$ or $B_{0}$, but not both, involved $\omega_{p}$ covers at
most $a+b$ patterns of Type 1--9.

\item[(2)] With both $A_{0}$ and $B_{0}$ involved $\omega_{p}$ covers at most $1$
pattern of Type $2-9$.

\item[(3)] With $A_{1}$ or $B_{1}$ (or both), but not $A_{0}$ or $B_{0}$, involved
$\omega_{p}$ covers at most $a+b-1$ patterns of Type  7--9.
\end{itemize}

Letting now $d_{i,j}$ denote the number of heads of Type $j$ patterns covered
by $\omega_{p}$ in the way described in category $(i)$ where $1\leq
i\leq3,1\leq j\leq9$ and using (\ref{l1}) the bounds for the sums $\sum
_{j}d_{i,j}$ and the zero $d_{i,j}$'s we have, as in the proof of Proposition
6,
\begin{eqnarray*}
\left|  T_{p}\right|  -\gamma H_{p}&\leq&1+\lambda_{1}\sum_{j}d_{1,j}%
+\lambda_{2}\sum_{j}d_{2,j}+\lambda_{9}\sum_{j}d_{3,j}-\gamma(a+b+2) \\
&\leq&\hbox{ }(9-16\gamma)(a+b)-(5-9\gamma)\leq13-23\gamma<0
\end{eqnarray*}
since $a+b\geq2$. This completes the proof.
\enddemo

{\it Remark}. The above proofs explain why we have only considered only
those nine types of patterns. For example it is now easy to show that if
$\omega_{p}$ covers the head of a pattern looking like $E_{2,2}\rightarrow
(E_{1},E_{1},E_{1},\ast)$ (which has not been included) then $L_{p}$ and
$R_{p}$ will have the properties mentioned in the above propositions.

\section{Good pairs}

We will say that a pair $(A,B)$ of a right interval $A\in{\cal F}^{+}(\mu)$
and a left interval $B\in{\cal F}^{-}(\mu)$ is {\it good} if there
exists $\omega_{p}\subseteq A\cap B$ such that $A=L_{p}$, $B=R_{p}$ and
\begin{equation}
\left|  T_{p}\right|  -\gamma H_{p}>0.%
\end{equation}
Using Propositions 6 and 7 we now conclude that any good pair $(A,B)$ must
satisfy the following:
\begin{itemize}
\item[(i)] Both $A$ and $B$ contain places from $E_{1}$.

\item[(ii)] There exists $\omega_{s}\subseteq A\cap B$ such that $(\omega_{s},A,B) $
covers an $\omega_{t}\in E_{1}$.
\end{itemize}

Suppose now that $(A,B)$ is a good pair. Then clearly $A$ uniquely determines
$B$ and vice versa. We define
\begin{equation}
w(A)=\min(A\cap 
{\textstyle\bigcup}
E_{1})<w(B)=\max(B\cap 
{\textstyle\bigcup}
E_{1}). 
\end{equation}
Clearly by (i) above we must have $\omega_{p}\subseteq A\cap B\subseteq
(w(A),w(B))$. Moreover we have the following.

\specialnumber{8}\proclaim{Lemma}
Suppose $(A,B)$ is a good pair. Then\/{\rm :}
\begin{itemize}
\item[{\rm (i)}] No $\omega_{q}\subseteq\lbrack\frak{l}(A),w(A)]\cup\lbrack w(B),\frak{r}%
(B)]$ can be the head of a Type $j$ pattern for any $1\leq j\leq9$.

\item[{\rm (ii)}] For every $\omega_{q}\subseteq\lbrack w(A),w(B)]$ we have $%
{\textstyle\bigcup}
T_{q}\subseteq\lbrack w(A),w(B)]$.

\item[{\rm (iii)}] Suppose that $\omega_{q}\subseteq\lbrack\frak{l}(A),w(A)]$ covers the
head of a Type $j$ pattern for some $1\leq j\leq9$. Then this can happen only
through the involvement of $L_{q}${\rm ,} which is then uniquely determined. A
symmetrical statement holds if $\omega_{q}\subseteq\lbrack w(B),\frak{r}(B)]$.
{\rm (}\/Here $R_{q}$ must be involved\/{\rm .)}
\end{itemize}

\endproclaim

\demo{Proof}
(i) Suppose $\omega_{q}\subseteq\lbrack\frak{l}(A),w(A)]$. Clearly $h_{p}%
\geq2$ and $\omega_{q}\subseteq A$. Using Lemma 6 it easily follows that there
must exist a left interval $I_{1}$ such that $(\omega_{q},A,I_{1})$ covers an
$\omega_{q_{1}}$ that is the head of a Type $1,2,3$ or $6$ pattern. Since
$q\leq w(A)$ and $[w(A),w(A)+1]\in E_{1}$ we must have $\frak{r}(I_{1})\leq
w(A)$ and therefore $\omega_{q_{1}}\subseteq\lbrack\frak{l}(A),w(A)]$. Arguing
similarly there must exist an $\omega_{q_{2}}\subseteq\lbrack\frak{l}%
(A),w(A)]$ (covered by $\omega_{q_{1}}$) that is the head of a Type $1$ or $2$
pattern and hence an $\omega_{q_{3}}\subseteq\lbrack\frak{l}(A),w(A)]\cap
\bigcup E_{1}$, which is a contradiction. The proof for $[w(B),\frak{r}(B)]$
is similar.
\vglue6pt 
(ii) Let $G,H$ be a pair such that $(\omega_{q},G,H)$ covers $\omega_{s}$
which is the head of some pattern. It is clear that $G\cup H\subseteq A\cup B
$ and so $\omega_{s}\in\lbrack\frak{l}(A),\frak{r}(B)]$. But also by (i)
$\omega_{s}$ cannot be contained in $[\frak{l}(A),w(A)]\cup\lbrack
w(B),\frak{r}(B)]$. Hence $\omega_{s}\subseteq\lbrack w(A),w(B)]$ and this
completes the proof.
\vglue6pt 
(iii) Suppose that there is a right interval $I$ different from $L_{q}$ and so
with $\frak{l}(L_{q})<\frak{l}(I)$ and a left interval $H$ such that
$(\omega_{q},I,H)$ covers an $\omega_{s}$ which is the head of some pattern.
As in (i) $\frak{r}(H)\leq w(A)$ and so $\omega_{s}\subseteq\lbrack
\frak{l}(I),w(A)]$. However (i) now implies that $\omega_{s}\subseteq
\lbrack\frak{l}(I),\frak{l}(A)]\subseteq I$. As in (i) there must exist a
right interval $H_{1}$ such that $(\omega_{s},I,H_{1})$ covers an
$\omega_{s_{1}}$ which is the head of a Type $1,2,3$ or $6$ pattern. Again we
get $\omega_{s_{1}}\subseteq\lbrack\frak{l}(I),w(A)]$ and so by (i)
$\omega_{s_{1}}\subseteq\lbrack\frak{l}(I),\frak{l}(A)]\subseteq I$. Now as in
(i) there must exist an $\omega_{t}\subseteq\lbrack\frak{l}(I),\frak{l}%
(A)]\cap\bigcup E_{1}$ and this is a contradiction since $\chi_{I}+\chi
_{L_{q}}=2$ on $[\frak{l}(I),\frak{l}(A)]$. Thus in any such covering $L_{q}$
must be involved$.$\pagebreak

To show that $L_{q}$ is uniquely defined suppose that for some other
$\omega_{q^{\prime}}\subseteq\lbrack\frak{l}(A),w(A)]$ that covers the head of
some pattern we had $L_{q}\neq L_{q^{\prime}}$. We may assume that
$\frak{l}(L_{q})<\frak{l}(L_{q^{\prime}})$. Then as before $\omega_{q^{\prime
}}$ must cover the head $\omega_{s}$ of some pattern, where $\omega
_{s}\subseteq$ $[\frak{l}(L_{q^{\prime}}),\frak{l}(A)]$ and this leads to a
similar contradiction. Hence $L_{q}$, if it exists, is uniquely defined.
\enddemo

{\it Remark}. If an $\omega_{q}$ as in Lemma 8(iii) exists then it is easy
to see that there is no left interval $G$ such that $(L_{q},G)$ is a good
pair. Indeed if such a $G$ existed then $L_{q}\cap G\subseteq\lbrack
w(L_{q}),w(G)]$ and so since $\omega_{q}\in L_{q}\cap A$ we must have
$G=J_{r}^{-}$ for some $r$ with $a_{r}<w(A)$ which implies that $G\subseteq
L_{q}\cap A$ and this is a contradiction.
\vglue8pt

Suppose now that $(A,B)$ is a good pair and define
\begin{equation}
T(A,B)=\{\omega_{s}:\omega_{s}\subseteq A\cup B\}
\end{equation}
and so $\left|  T(A,B)\right|  =\frak{r}(B)-\frak{l}(A)=\left|  A\cup
B\right|  $.

Next we consider $A$. If $A$ is clean then let $g(A)=K(A)=K^{\ast}(A)=0$. If
$A$ is not clean then we write (see \S 3) $P\cap(\frak{l}%
(A),w(A)]=\{a_{s},\ldots , a_{t}\}\neq\emptyset$ (and so $A=J_{s-1}^{+}$) and
with
\begin{equation}
g(A)=a_{t}-\frak{l}(A)\hbox{ and }K(A)=k_{s}+\cdots + k_{t}\label{nc1}%
\end{equation}
we now define $K^{\ast}(A)$ as follows:
\vglue4pt
(i) if there exists at least one $\omega_{q}\subseteq\lbrack\frak{l}(A),w(A)]$
as in the statement of Lemma 8(iii), $K^{\ast}(A)$ is equal to the total
number of bricks that correspond to the left intervals $J_{t}^{-}%
,\ldots , J_{s}^{-}$ or to right intervals $J_{l}^{+}$ with $l<i$ and lie over
$[\frak{l}(A),w(A)]\backslash L_{q}$ plus the length of the interval
$L_{q}\cap A$ (note that $L_{q}$ is uniquely determined and that we must have
$\frak{r}(L_{q})\leq w(A)$), and
\vglue4pt
(ii) if no such $\omega_{q}$ exists, $K^{\ast}(A)$ is equal to the total
number of bricks that lie over $[\frak{l}(A),w(A)]$ and correspond to either
the left intervals $J_{t}^{-},\ldots , J_{s}^{-}$ or to right intervals $J_{l}^{+}
$ with $l<i$.
\vglue4pt
Note that in both cases bricks that correspond to $A$ are not counted in
$K^{\ast}(A)$.

We also consider $B$ and define $g(B),K(B),K^{\ast}(B)$ in a completely
symmetrical way.

Regarding the masses that lie in $(w(A),w(B))$ we set
\begin{equation}
K(A,B)=\sum_{w(A)<a_{r}<w(B)}k_{r}%
\end{equation}
and now we define
\begin{equation}
H(A,B)=\left|  A\right|  +K^{\ast}(A)+K(A)+2K(A,B)+K(B)+K^{\ast}(B)+\left|
B\right|  .\quad
\end{equation}

It is easy to see that by our construction
\begin{equation}
H(A,B)\leq\sum_{\omega_{p}\subseteq A\cup B}h_{p}\label{h}.
\end{equation}
(For example if $a_{r}\in(w(A),w(B))$ then we must have $J_{r}\subseteq
(w(A),w(B))$ and so all the $2k_{r}$ bricks corresponding to $J_{r}$ lie over
$A\cup B$.) Also if $A$ is {\it not clean} then $K^{\ast}(A)>0$ and each
place in $[\frak{l}(A),w(A)]$ contributes at least {\it two} bricks in
$H(A,B)$ (one from $A$ and at least one counted in $K^{\ast}(A)+K(A)$), in
particular $K^{\ast}(A)+K(A)\geq g(A)$.

The main thing now is to prove the following basic.

\specialnumber{8}\proclaim{Proposition}
There exists at least one good pair $(A,B)$ such that
\begin{equation}
\left|  T(A,B)\right|  >\gamma H(A,B).\label{g5}%
\end{equation}

\demo{Proof}
First of all we have the following.

\specialnumber{}\proclaim{Lemma}
Given any two good pairs $(A,B)$ and $(A^{\prime},B^{\prime})$ with
$\frak{l}(A)<\frak{l}(A^{\prime})$ we must have
\begin{equation}
\frak{r}(B)\leq\frak{l}(A^{\prime}).%
\end{equation}
\endproclaim

\demo{Proof}
Assume that $A=J_{i}^{+}$, $B=J_{j}^{-}$, $A^{\prime}=J_{s}^{+}$ and
$B^{\prime}=J_{r}^{-}$ and moreover that $a_{i}<a_{s}$ but $a_{j}>a_{s}$. We
must have $a_{j}<a_{r}$; otherwise, since both $A\cap B$ and $A^{\prime}\cap
B^{\prime}$ are nonempty we would have $\chi_{A}+\chi_{B}+\chi_{A^{\prime}%
}+\chi_{B^{\prime}}\geq2$ on $[a_{s},a_{r}]$ contradiction. Considering now
the symmetric of $A^{\prime}$ and $B$ intervals $H=J_{s}^{-} $ and
$G=J_{j}^{+}$ we have $\chi_{H}+\chi_{A^{\prime}}+\chi_{B^{\prime}}\geq1$ on
$[a_{s}-k_{s},a_{r}]$ and $\chi_{G}+\chi_{A}+\chi_{B}\geq1$ on $[a_{i}%
,a_{j}+k_{j}]$. Consider now an $\omega_{q}\in E_{1}$ contained in $B$. Then
we must have $q\leq a_{s}-k_{s}$ and so $a_{j}-k_{j}=\frak{l}(B)<q\leq
a_{s}-k_{s}$. In a similar way we obtain $a_{s}+k_{s}>a_{j}+k_{j}$. These give
$a_{j}-a_{s}<k_{j}-k_{s}<a_{s}-a_{j}$ contradiction since $a_{j}>a_{s}$. This
completes the proof.
\enddemo

In view of the above lemma we can number all the good pairs of $\mu$ (if any)
as $(A_{1},B_{1}),\ldots , (A_{d},B_{d})$ so that $\frak{r}(B_{i})\leq
\frak{l}(A_{i+1})$ for $i=1,\ldots , d-1$. This implies that the sets
$T(A_{1},B_{1}),\ldots , T(A_{d},B_{d})$ are pairwise disjoint. Let
\begin{equation}
W=\{\omega_{1},\ldots , \omega_{N}\}\backslash\bigcup_{i=1}^{d}T(A_{i},B_{i})
\end{equation}
and consider the collection ${\cal S}$ of all $T_{p}$'s where either: (i)
$\omega_{p}\in W$ and is not the head of any Type $j$ pattern for any $1\leq
j\leq9$ or (ii) $\omega_{p}$ is the head of some such pattern but there is
$1\leq i\leq d$ such that $\omega_{p}$ is covered by an $\omega_{q}%
\subseteq\lbrack\frak{l}(A_{i}),w(A_{i})]\cup\lbrack w(B_{i}),\frak{r}%
(B_{i})]$ (through the involvement of $L_{q}$). We then have the following.

\specialnumber{10}\proclaim{Lemma}
{\rm (i)} Any $T_{p}\in{\cal S}$ is disjoint from $\bigcup_{i=1}^{d}T(A_{i}%
,B_{i})$.
\vglue4pt
{\rm (ii)} We have
\begin{equation}
E_{1}\subseteq\bigcup_{T_{p}\in{\cal S}}T_{p}\cup\bigcup_{i=1}^{d}%
T(A_{i},B_{i}).%
\end{equation}

{\rm (iii)} For every $T_{p}\in{\cal S}$ we have $\left|  T_{p}\right|  <\gamma
H_{p}$.
\endproclaim

\demo{Proof}
(i) Suppose that $T_{p}\in{\cal S}$ and $\omega_{q}\in T_{p}\cap
T(A_{i},B_{i})$ for some $i$. If $q=p$ then Lemma 8 and the definition of
${\cal S}$ easily imply that $T_{p}\notin{\cal S}$. If $q\neq p$ then
$\omega_{q}$ is the head of some Type $j$ pattern and so by Lemma 8 we must
have $\omega_{q}\subseteq\lbrack w(A_{i}),w(B_{i})]$. But then it is easy to
see that $\omega_{q}$ can be covered only if $A_{i}$,$B_{i}$ or some of the
masses corresponding to positions in $[w(A_{i}),w(B_{i})]$ are involved and
this would give $\omega_{c(q)}\in T(A_{i},B_{i})$. Continuing this (for at
most three steps) we conclude that $\omega_{p}\in T(A_{i},B_{i})$ which as we
have seen is a contradiction.
\vglue4pt
(ii) Suppose that $\omega_{q_{0}}\in E_{1}\backslash\bigcup_{i=1}^{d}%
T(A_{i},B_{i})$ and let $q_{1}=c(q_{0}),q_{2}=c(q_{1}),\ldots$ (that is
$\omega_{q}$ is covered by $\omega_{q_{1}}$ which is covered by $\omega
_{q_{2}}$ and so on). Clearly $h_{q_{r}}\geq2$ for all $r\geq1$. Let $m\geq1$
be the smallest possible integer such that $\omega_{q_{m}}$ is not the head of
a Type $j$ pattern for any $1\leq j\leq9$ (note that $\omega_{q_{0}}$ is the
head of a Type $1$ pattern). Such an $m$ exists since each such pattern
contains at most eight places and by Lemma 6 no cycles (that is chains of the
form $\omega_{p_{1}}\rightarrow\omega_{p_{2}}\rightarrow\cdots\rightarrow
\omega_{p_{s}}=\omega_{p_{1}}$). By Lemma 8 we conclude that $\omega_{q_{r}%
}\in W$ for all $0\leq r\leq m-1$. If $\omega_{q_{m}}\in T(A_{i},B_{i})$ for
some $i$ then we must have $\omega_{q_{m}}\subseteq\lbrack\frak{l}%
(A_{i}),w(A_{i})]\cup\lbrack w(B_{i}),\frak{r}(B_{i})]$  (otherwise Lemma
8(ii) would imply that $\omega_{q_{0}}\in T_{q_{m}}\subseteq T(A_{i},B_{i})$)
and so $\omega_{q_{0}}\in T_{q_{m-1}}\in{\cal S}$. If $\omega_{q_{m}}\in W$
then $\omega_{q_{0}}\in T_{q_{m}}\in{\cal S}$.
\vglue4pt
(iii) Consider $T_{p}\in{\cal S}$ . Suppose that $\omega_{p}\in W$ is not
the head of any Type $j$ pattern. Then by (i), $(L_{p},R_{p})$ is not a good
pair hence we have $\left|  T_{p}\right|  <\gamma H_{p}$. If $\omega_{p}\in W$
is the head of such a pattern then the definition of $H_{p}$ (see (\ref{h1a}))
shows that $\left|  T_{p}\right|  <\gamma H_{p}$.
\enddemo

We next let
\begin{equation}
D=\{\omega_{1},\ldots , \omega_{N}\}\backslash\left(  \bigcup_{T_{p}\in{\cal S}%
}T_{p}\cup\bigcup_{i=1}^{d}T(A_{i},B_{i})\right)
\end{equation}
and note that by Lemma 10(ii) we have $h_{q}\geq2$ for every $\omega_{q}\in
D$. Then by letting $T_{p_{1}},\ldots , T_{p_{m}}$ be all the maximal $T_{p}$'s
from ${\cal S}$, which by Lemma 7 are pairwise disjoint and cover
$\bigcup_{T_{p}\in{\cal S}}T_{p}$ we have
\begin{equation}
\left|  J(\mu)\right|  =N=\sum_{r=1}^{m}\left|  T_{p_{r}}\right|  +\sum
_{i=1}^{d}\left|  T(A_{i},B_{i})\right|  +\left|  D\right|  .\label{g1}%
\end{equation}

Now the following holds.

\specialnumber{11}\proclaim{Lemma}
We have
\begin{equation}
\sum_{p=1}^{N}h_{p}\geq\sum_{r=1}^{m}H_{p_{r}}+\sum_{i=1}^{d}H(A_{i}%
,B_{i})+2\left|  D\right|  .\label{g2}%
\end{equation}
\endproclaim

\demo{Proof}
It is enough to show that the right-hand side of (\ref{g2}) is at most as
large as the total number of bricks that lie over all $\omega_{s}$'s. Using
that $h_{p}\geq2$ for all $\omega_{p}\in D$, Lemma 9, Lemma 10(i), (\ref{h}),
the remark following Lemma 8 and the definitions of the $H_{p}$'s and the
$H(A_{i},B_{i})$'s we easily  see that the only case that should be
considered is when $\omega_{p}$ is the head of a Type $j$ pattern and is
covered by an $\omega_{q}\subseteq\lbrack\frak{l}(A_{i}),w(A_{i})]\cup\lbrack
w(B_{i}),\frak{r}(B_{i})]$ for some $i$ in which case $H_{p}$ counts one more
brick than the ones involved. Assume $\omega_{q}\subseteq\lbrack\frak{l}%
(A_{i}),w(A_{i})]$. Then $L_{q}$ is uniquely determined and $\omega_{q}$ can
cover at most as many such heads $\omega_{p}$ as there are bricks lying over
$\omega_{q}\subseteq L_{q}$ that correspond to {\it left} intervals whose
right endpoints are contained in $[\frak{l}(A_{i}),w(A_{i})]$. However by the
definition of $K^{\ast}(A_{i})$ it is clear that all these bricks are not
counted in $H(A_{i},B_{i})$. A similar reasoning for the case $\omega
_{q}\subseteq\lbrack w(B_{i}),\frak{r}(B_{i})]$ completes the proof of
(\ref{g2}).
\enddemo

Now since $R(\mu)>1+\gamma$ we have $2\gamma\sum_{p=1}^{N}h_{p}\leq2\gamma
K_{1}^{n}<\left|  J(\mu)\right|  $ and so using Lemma 10(iii), (\ref{g1}) and
(\ref{g2}) we conclude that there must exist at least one $i$ (hence at least
one good pair) such that $\left|  T(A_{i},B_{i})\right|  >\gamma H(A_{i}%
,B_{i})$. This completes the proof of the proposition. \hfill\qed
\enddemo

\section{The core of a good pair}

Now, using the theorem, we can find and fix a good pair $(A,B)$ that satisfies
(\ref{g5}).

\specialnumber{12}\proclaim{Lemma}
The interval $A\cap B$ {\rm (}\/corresponding to the pair $(A,B)$\/{\rm )} cannot cover places
in both $A\cap\bigcup E_{1}$ and $B\cap\bigcup E_{1}$. Moreover if it covers
at least one place in $A\cap\bigcup E_{1}$ then it cannot cover any place in
$B\backslash A$.
\endproclaim

\demo{Proof}
Suppose $A=J_{i}^{+}$ and $B=J_{j}^{-}$. Then clearly $I_{i,j}(\mu)$ is a
special interval; therefore $\left|  I_{i,j}(\mu)\right|  =\left|  A\cap
B\right|  $ and so $A\cap B$ is placed, without breaking it, over $E(\mu)$.
Going to the gap interval $J(\mu)$ if $x,y\in J(\mu)$ are covered by $A\cap B$
then since $Q$ is distance nondecreasing we must have $\left|  x-y\right|
\leq\left|  A\cap B\right|  $. However if $\omega_{p}\subseteq A\cap\bigcup
E_{1}$ and $\omega_{q}\subseteq B\cap\bigcup E_{1}$, or $\omega_{q}\subseteq
B\backslash A$, then it is easy to see that $\left|  q-p\right|  >\left|
A\cap B\right|  $ and this completes the proof.
\enddemo

In view of the above lemma and the properties shared by any good pair we may
assume that $A\cap B$ covers at least one place in $A\cap\bigcup E_{1}$ and so
no place in $B\cap\bigcup E_{1}$ or $B\backslash A$.

We then let
\begin{equation}
z=z(A,B)=\max\{s:\omega_{s}\subseteq A\cap\bigcup E_{1}\}>w(A).%
\end{equation}
Since $A\cap B$ does not contain any place from $E_{1}$ we have $A\cap
B\subseteq\lbrack z(A,B),w(B)]$. Next we write
\begin{equation}
\lbrack z(A,B),w(B)]\cap P=\{a_{p},a_{p+1},\ldots , a_{q}\}
\end{equation}
and define the {\it core} of $(A,B)$ to be the measure
\begin{equation}
\sigma=\sigma(A,B)=\sum_{r=p}^{q}k_{r}\delta_{y_{r}}%
\end{equation}
that corresponds to these positions.

\demo{{R}emark} (i) The set $[z(A,B),w(B)]\cap P$ must be nonempty. If it
were empty then $B$ would not  interact with any right interval other than $A$
to the left of $w(B)$ and also $A$ would not interact with any left interval
other than $B$ after $z$. This would imply that $A\cap B\subseteq\lbrack
z(A,B),w(B)]$ must be covered only by the intersection $A\cap B$. But this is
impossible since $A\cap B$ must cover at least one place in $E_{1}$ and this
place must be outside $A\cap B$.
\vglue4pt
(ii) Note the nonsymmetrical way with respect to $A$ and $B$ the core interval
is defined (a $\max$ for right intervals would correspond to a $\min$ for left
intervals). This is forced because of the location of the special interval
corresponding to $(A,B)$ (see also the construction in the Appendix).
\enddemo

We will now show that without affecting the core of $(A,B)$ we may assume that
both intervals $A$ and $B$ are {\it clean}. This would be important in the
next section and is furnished by the following.

\specialnumber{9}\proclaim{Proposition}
For the good pair $(A,B)$ considered above there exists an admissible measure
$\bar{\mu}$ {\rm (}\/which in general might contain more positions
than~$\mu$\/{\rm )} and a good pair $(\bar{A},\bar{B})$ associated to the families
${\cal F}^{\pm }(\bar{\mu})$ corresponding to the gap interval of $\bar{\mu}$ such
that\/{\rm :}
\begin{itemize}
\ritem{(i)} $\left|  T(\bar{A},\bar{B})\right|  >\gamma H(\bar{A},\bar{B})$.

\ritem{(ii)} Both the right interval $\bar{A}$ and the left interval $\bar{B}$ are clean.

\ritem{(iii)} The core $\sigma(\bar{A},\bar{B})$ of the good pair $(\bar{A},\bar{B}) $
is identical to the core $\sigma(A,B)$ of $(A,B)$.

\ritem{(iv)} For any measure $\nu$ formed from masses of $\bar{\mu}$ whose associated
positions in $J(\bar{\mu})$ are contained in the interior of $\bar{A}\cup
\bar{B}$ we have $\left|  E(\nu)\right|  \leq2(1+\gamma)\left\|  \nu\right\|
$.
\end{itemize}

\endproclaim

\demo{Proof}
If both $A$ and $B$ are clean there is nothing to prove. Suppose that $A$ is
not clean. Define then $w(A),g=g(A),K=K(A)$ and $K^{\ast}=K^{\ast}(A)$ as in
Section 5, write $A=J_{p}^{+}$ and suppose that for some $i>p$
\begin{equation}
(w(A),+\infty)\cap P=\{a_{i+1},\ldots , a_{n}\}
\end{equation}
(it is obviously nonempty) and so $K=K(A)=K_{p+1}^{i}$. We will not change
anything in the part of the gap interval of $\mu$ that lies to the right of
$w(A)$. Let $s\geq i$ be such that
\begin{equation}
a_{s}\leq\frak{l}(A)=a_{p}+k_{p}<a_{s+1}.
\end{equation}
Then the considerations in Section 3 and Lemma 4 imply that $[y_{i}%
+k_{i},y_{i}+k_{i}+1]$ is covered by a special interval $I_{p,t}(\mu)$ for
some $t>s$ and moreover using the remark following Lemma 4 we may and will
assume that
\begin{equation}
y_{i}+k_{i}=\frak{l}(I_{p,t}(\mu))=y_{t}-k_{p}-K-K_{i+1}^{t}%
\end{equation}
and so
\begin{equation}
g+K=\mathop{\rm dist}\nolimits (A,a_{t})+K_{i+1}^{t-1}.\label{p1}%
\end{equation}
Now we fix an admissible measure $\tau$ all whose entries are rational numbers
such that
\begin{eqnarray}
E(\tau)&=&[y_{i}+k_{i}-\mathop{\rm dist}\nolimits (A,a_{t})-K_{s+1}^{t-1},y_{i}%
+k_{i}],\label{p2} \\
\left|  E(\tau)\right| & =&2(1+\gamma-\varepsilon)\left\|  \tau\right\|
,%
\end{eqnarray}
where $\varepsilon>0$ is small to be fixed later and such that the maximum
(individual) mass appearing in the positions of $\tau$ is so small that no
mass of $\tau$ interacts with any $k_{r}\delta_{r}$ for any $r>i$. Such a
measure can be constructed for example by the proceedure that leads to the
lower bound for $C$ (see \cite{Mel} or the Appendix here) and an appropriate scaling-translation.

Let
\begin{equation}
\bar{K}=\left\|  \tau\right\|  \hbox{ and }\bar{g}=\left|  G(\tau)\right|
=\left|  E(\tau)\right|  -2\bar{K}=2(\gamma-\varepsilon)\bar{K}=\left|
J(\tau)\right|  .\hskip.25in
\end{equation}
Next we define
\begin{equation}
\bar{k}_{p}=k_{p}+K-\bar{K},%
\end{equation}
noticing that $\bar{k}_{p}>0$ since $2\bar{K}<g+K<k_{p}+K$.

Consider now the measure
\begin{equation}
\bar{\mu}=\bar{k}_{p}\delta_{y_{p}}+\tau+\sum_{r=i}^{n}k_{r}\delta_{y_{r}%
}.%
\end{equation}
Here the index $p$ is used for convenience only, since we have no control on
the number of positions in $\tau$. Consequently we will not associate indices
to the positions of $\tau$.

Also by multiplying all entries in $\mu$ and $\bar{\mu}$ by the same
appropriately chosen large integer we may assume that all such entries are
{\it integers}.

Also consider in the gap interval $J(\bar{\mu})$ the pair $(\bar{A},B)$ ($B$
as before) where
\begin{equation}
\bar{A}=J_{p}^{+}(\bar{\mu})
\end{equation}
is the right interval corresponding to $\bar{k}_{p}\delta_{y_{p}}$. We will
show that $\bar{\mu}$ is admissible, that the pair $(\bar{A},B)$ is good with
$\bar{A}$ clean and also that (i), (iii) and (iv) are satisfied. This will
actually complete the proof since in case $B$ is also not clean we can apply a
similar symmetrical construction with $B$ and the measure $\bar{\mu}$ to
satisfy all conditions.

Since $y_{p}$ and $y_{t}$ have not been altered and since $\bar{k}_{p}+\bar
{K}=k_{p}+K(A)$ we have
\begin{equation}
I_{p,t}(\mu)=I_{p,t}(\bar{\mu}).%
\end{equation}
Consequently in view of Lemma 2 and since $I_{p,t}(\mu)$ is a special interval
we conclude that $I_{p,t}(\bar{\mu})$ must also be a special interval (with
respect to $\bar{\mu}$) and therefore in the gap intervals $J(\mu)$ and
$J(\bar{\mu})$ the right endpoints $\frak{r}(A)$ and $\frak{r}(\bar{A})$ must
respectively be located at the same point of $J_{t}^{-}(\mu)$ and $J_{t}%
^{-}(\bar{\mu})$. This in view of Proposition 2(ii) and Lemma 3 and, since we
have not altered $\mu$ to the right of $y_{i}$, implies that we must have
\begin{equation}
I_{p,r}(\mu)=I_{p,r}(\bar{\mu})
\end{equation}
for every $r\geq t$ and since (the nonempty of) these intervals together with
$E(\sum_{r=i+1}^{n}k_{r}\delta_{r})$ cover the space $[y_{i}+k_{i},y_{n}]$ of
$E(\mu)$ (note that $I_{l,r}(\mu)=\emptyset$ if $l\leq i<r$ with $l\neq p$
and that the nonempty, if any, of the intervals $I_{p,r}(\mu)$ for $r<t$ are
located to the right of $y_{i}+k_{i}$) we conclude that
\begin{equation}
\lbrack y_{i}+k_{i},y_{n}]\subseteq E(\bar{\mu}).%
\end{equation}
Also it is clear that $E(\tau)\subseteq E(\bar{\mu})$. Now as remarked above
in the gap interval of $\bar{\mu}$ the interval $\bar{A}=J_{p}^{+}(\bar{\mu})$
must contain all positions that correspond to the masses $k_{i+1}\delta
_{i+1},\ldots , k_{s}\delta_{s}$ (and obviously all the positions corresponding to
$\tau$) we have
\begin{equation}
\frak{r}(F_{p}(\bar{\mu}))=y_{p}+\bar{k}_{p}+\bar{K}+K_{i+1}^{s}=y_{p}%
+k_{p}+K+K_{i+1}^{s}. \hskip.5in
\end{equation}
Hence in view of (\ref{p1}) and (\ref{p2})
\begin{equation}
\frak{l}(E(\tau))-\frak{r}(F_{p}(\bar{\mu}))=y_{i}+k_{i}-g-K-y_{p}-k_{p}-K=0\hskip.5in
\end{equation}
and this now implies that $E(\bar{\mu})$ is connected, therefore that
$\bar{\mu}$ is admissible (the separability inequalities being here obvious).

Now by the way $\tau$ is chosen (iv) is satisfied and also, since nothing has
changed after $w(A)$, it is clear, using also Lemma 8(ii), that the pair
$(\bar{A},B)$ is good and that its core satisfies $\sigma(\bar{A}%
,B)=\sigma(A,B)$.

To prove (i) we form the gap intervals of $\mu$ and $\bar{\mu}$ simultaneously
shrinking the corresponding central intervals $I_{r,r}$ of $\mu$ and $\bar
{\mu}$\ in such a way that in both cases the point $b=y_{i}+k_{i}$ is kept
fixed. In this way in both gap intervals the segments that lie in
$[b,+\infty)$ are identical and also $\frak{r}(A)$ $=\frak{r}(\bar{A})$. Now
in $\mu$, as we already know, a gap of exactly $g$ will be formed between
$\frak{l}(A)$ and $b$. In $\bar{\mu}$ however $E(\tau)$ will shrink to the
interval $[b-\left|  J(\tau)\right|  ,b]$ and between $\frak{l}(\bar{A})$ and
$b-\left|  J(\tau)\right|  $ a gap of exactly $\bar{K}+K_{i+1}^{s}$ will be
formed, proving thus that in particular $\bar{A}$ is clean (since the
individual masses of $\tau$ have been chosen very small). Hence it is easy to
see that
\begin{equation}
X=\left|  T(\bar{A},B)\right|  -\left|  T(A,B)\right|  =\bar{K}+K_{i+1}%
^{s}+\bar{g}-g
\end{equation}
and
\begin{equation}
Y=\left|  H(\bar{A},B)\right|  -\left|  H(A,B)\right|  =\bar{K}+K_{i+1}%
^{s}+\bar{g}+2\bar{K}-g-(K+K^{\ast}).\qquad
\end{equation}
In view of (\ref{g5}) to prove (i) it is enough to show that $X>\gamma Y$. We
have
\begin{equation}
X-\gamma Y=(1-\gamma)K_{i+1}^{s}+\gamma(K+K^{\ast})-(1-\gamma)g+(1-\gamma
)\bar{g}-(3\gamma-1)\bar{K}.\label{p3}%
\end{equation}
Using (\ref{p2}) it is now easy to compute  that
\begin{equation}
(1-\gamma)\bar{g}-(3\gamma-1)\bar{K}=\left(\frac{1}{2}-\gamma-\varepsilon^{\prime
}\right)(K_{s+1}^{t-1}+\mathop{\rm dist}\nolimits (A,a_{t}))\hskip.5in \label{p4}%
\end{equation}
where $\varepsilon^{\prime}=\dfrac{\varepsilon}{2(\gamma-\varepsilon+1)} $.
Moreover we have $K+K^{\ast}\geq g$, since obviously each place in $g$
contributes at least one brick counted in $K+K^{\ast}$. Hence (since
$\gamma>\frac{1}{2}$),
\begin{equation}
\gamma(K+K^{\ast})-(1-\gamma)g\geq\left(\gamma-\frac{1}{2}\right)(g+K+K^{\ast}%
).\label{p5}%
\end{equation}
Now using  (\ref{p4}), (\ref{p5}) and (\ref{p1}) in (\ref{p3}) and observing
that we must have $K^{\ast}>0$ we get
\begin{eqnarray} \qquad
X-\gamma Y&\geq&(1-\gamma)K_{i+1}^{s}+\left(\frac{1}{2}-\gamma-\varepsilon^{\prime
}\right)(K_{s+1}^{t-1}+\mathop{\rm dist}\nolimits (A,a_{t})) \\[5pt]
&&+\ \left(\gamma-\frac{1}{2}\right)(K_{i+1}^{t-1}+\mathop{\rm dist}\nolimits (A,a_{t})+K^{\ast
}) \nonumber\\[5pt]
&=&\frac{1}{2}K_{i+1}^{s}+\left(\gamma-\frac{1}{2}\right)K^{\ast}-\varepsilon^{\prime
}(K_{s+1}^{t-1}+\mathop{\rm dist}\nolimits (A,a_{t}))>0  \nonumber
\end{eqnarray}
if $\varepsilon>0$ has been choosen small enough. This \pagebreak completes the
proof.
\enddemo

\section{The basic estimate for the core}

We will now consider a good pair $(A,B)$ in which both $A$ and $B$ are clean
and is such that (\ref{g5}) is satisfied. This pair can be a part of $\mu$ or
be produced as in Proposition 9. In both cases its core $\sigma(A,B)$ is a
part of $\mu$ and contains less $n$ positions. For convenience we will change
the numbering of the $y_{i},a_{i}$ and $k_{i}$'s, introducing negative indices
and also introduce if necessary (at most) two positions in $\mu$ (or $\bar
{\mu}$) with masses $0$ in such a way that
\begin{equation}
\sigma=\sigma(A,B)=\sum_{i=1}^{m}k_{i}\delta_{y_{i}}\label{z1},
\end{equation}
where $m\leq n$ and moreover so that there are $1\leq r<s\leq n$ \ ($r<s$
since $A\cap B$ covers at least one place) with
\begin{equation}
\frak{r}(A)=a_{s}\hbox{ and }\frak{l}(B)=a_{r}.\label{z2}%
\end{equation}
It is easy to see that these new zero mass positions will not affect any of
the covering properties of ${\cal F}(\mu)$ or related estimates, but will
make our computations easier.

We will also use the following notation: For any $i<j$ we will let
\begin{equation}
\alpha_{i}^{j}=a_{j}-a_{i}\label{z3}%
\end{equation}
and we will let $\alpha_{i}^{j}=0$ if $j\leq i$.

Now the gap interval of $\sigma$ is $J(\sigma)=[a_{1},a_{m}]$. Doing that we
would have
\begin{equation}
A=J_{-p}^{+}\label{z4}%
\end{equation}
for some integer $p>0$ and we will also consider the intermediate measure
\begin{equation}
\nu=\sum_{i=-p+1}^{0}k_{i}\delta_{y_{i}}.\label{z5}%
\end{equation}

As for $B$ since it is also clean it is easy to see that Proposition 5 implies
that
\begin{equation}
B=J_{m+1}^{-}\hbox{ and }S=a_{m+1}-a_{m}\leq K_{r}^{m}.\label{z6}%
\end{equation}

We will now analyse $A$. Let
\begin{equation}
\rho=\left|  E(\nu)\backslash\bigcup_{i=-p+1}^{0}I_{i,i}(\nu)\right|  \hbox{
and }K=\sum_{i=-p+1}^{0}k_{i}=\left\|  \nu\right\|  .\label{z7}%
\end{equation}
Since $[a_{-p+1},a_{0}]\subseteq A$ is surrounded by places in $E_{1}$ we
conclude that no interval of ${\cal F}(\nu)$ interacts with any interval
other than $A$ and the interactions with $A$ produce an interval of length $K
$ in $F_{-p}\backslash(-\infty,y_{-p}+k_{-p}]$ (where $F_{-p}=F_{-p}(\mu)$).
Actually we have
\begin{equation}
\left|  F_{-p}\backslash(-\infty,y_{-p}+k_{-p}]\right|  =K+K_{1}^{s}%
.\label{z8}%
\end{equation}
The interval $I_{-p,m+1}(\mu)$ (or $\bar{\mu}$) that corresponds to $A\cap B $
can cover, by Lemma~12, only points $x\in G(\mu)$ such that $Q^{-1}(x)\in A$
and moreover it covers at least one place of $G(\mu)$ that corresponds to some
place in $A\cap\bigcup E_{1}$. In particular,
\begin{equation}
I_{-p,m+1}(\mu)\subseteq(-\infty,y_{s}].\label{z9}%
\end{equation}
Therefore denoting by $D$ the part of $I_{-p,m+1}(\mu)$ that lies in
$(-\infty,y_{1}]$ and also corresponds to the places in $(-\infty,a_{1})$ in
the gap interval covered by $A\cap B$ and by $h$ the part that lies in
$[y_{1},y_{s}]$, that is the, possibly empty, space in $[y_{1},y_{s}%
]\backslash E(\sigma)$ covered by $I_{-p,m+1}(\mu)$, we have (since $Q$ is
distance nondecreasing)
\begin{equation}
D>0\hbox{ and }D+h\leq\left|  A\cap B\right|  =\alpha_{r}^{s}.%
\label{z10}%
\end{equation}

Now we thus have $\frak{l}(I_{-p,m+1}(\mu))<y_{1}-k_{1}$ and by Lemma 3 we
see that 
\begin{equation}
g=\frak{l}(I_{-p,m+1}(\mu))-\frak{r}(F_{-p})=\mathop{\rm dist}\nolimits %
(a_{m+1},J_{-p})+K_{s+1}^{m}=\alpha_{s}^{m}+S+K_{s+1}^{m},\enspace\label{z11}%
\end{equation}
where (\ref{z11}) defines $g$. Hence by the above considerations and
Proposition 3(iii) the interval $(\frak{r}(F_{-p}),\frak{l}(I_{-p,m+1}(\mu)))$
in $E(\mu)$ must be covered by $E(\nu)$ and some of the nonempty special
intervals $I_{-p,j}(\mu)$ for $s+1\leq j\leq m$. Hence there is $\lambda\geq0$
such that
\begin{equation}
\lambda\leq\sum_{j=s+1}^{m}\left|  A\cap J_{j}^{-}\right|  \hbox{ and
}g=\lambda+\left|  E(\nu)\right|  =\lambda+\rho+2K.\label{z12}%
\end{equation}
This in turn implies that the total space in the gap interval $E(\mu)$,
between $a_{-p}+K+K_{1}^{s}$ and $a_{1}-D$, is at most $\lambda+\rho$. Hence
\begin{equation}
\alpha_{-p}^{1}=a_{1}-a_{-p}\leq K+K_{1}^{s}+D+\lambda+\rho.\label{z13}%
\end{equation}
Moreover since $\nu$ has less than $n$ positions (or see Proposition 9(iv)) we
have $\left|  E(\nu)\right|  \leq2(1+\gamma)N(\nu)$ and so
\begin{equation}
\dfrac{\rho}{2K}\leq\gamma.\label{z14}%
\end{equation}

Turning now to the core $\sigma$ we have that since no mass of $\sigma$
interacts with any mass outside $\sigma$ other than those corresponding to $A$
and $B$ and since all nonempty $I_{-p,j}(\mu)$ for $s+1\leq j\leq m$ are
situated to the left of $I_{-p,m+1}(\mu)$ whose left endpoint is smaller than
$y_{1}-k_{1}$, the interval $[y_{1}-k_{1},y_{m}+k_{m}]$ can be covered only by
$E(\sigma)$, the part $h$ of $I_{-p,m+1}(\mu)$ and possibly some of the
nonempty special intervals $I_{i,m+1}(\mu)$ for $1\leq i\leq r-1$. Hence
denoting by $u\geq0$ the measure of $[y_{1},y_{m}]\backslash(E(\sigma)\cup
I_{-p,m+1}(\mu))$ we have
\begin{equation}
u\leq\sum_{j=1}^{r-1}\left|  B\cap J_{j}^{+}\right|  \hbox{ and }(y_{m}%
+k_{m})-(y_{1}-k_{1})\leq\left|  E(\sigma)\right|  +u+h.\hskip.25in\label{z15}%
\end{equation}
Therefore since $\sigma$ contains less than $n$ (nonzero mass) positions we
have $\left|  E(\sigma)\right|  \leq2(1+\gamma)\left\|  \sigma\right\|  $ and
so
\begin{equation}
\alpha_{1}^{m}\leq2\gamma K_{1}^{m}+u+h.\label{z16}%
\end{equation}

Now to use the above information efficiently we introduce the estimate
(\ref{g5}) satisfied by the pair $(A,B)$. This gives
\begin{eqnarray}
a_{m+1}-a_{-p}  & =&T(A,B)>\gamma H(A,B) \label{z17} \\
& \geq&\gamma(a_{m+1}-a_{r}+a_{s}-a_{-p}+2K+2K_{1}^{m})\nonumber
\end{eqnarray}
and so
\begin{equation}
(1-\gamma)(S+\alpha_{1}^{m}+\alpha_{-p}^{1})>\gamma(2K+2K_{1}^{m}+\alpha
_{r}^{s}).\label{z18}%
\end{equation}
Using now the estimates (\ref{z10}) and (\ref{z13}) and since $\gamma<1$ we
get
\begin{eqnarray}
&&(1-\gamma)(K_{1}^{s}+S+\lambda+\alpha_{1}^{m}-h)+(1-\gamma)\rho-(3\gamma
-1)K\label{z19}\\
&&\qquad \quad >2\gamma K_{1}^{m}+(2\gamma-1)\alpha_{r}^{s}.\nonumber%
\end{eqnarray}

Moreover using (\ref{z12}) and (\ref{z14}) we may write
\begin{equation}
\rho=\eta(g-\lambda)\hbox{ and }K=\dfrac{1-\eta}{2}(g-\lambda),\label{z20}%
\end{equation}
where
\begin{equation}
\eta\leq\dfrac{\gamma}{\gamma+1},\label{z21}%
\end{equation}
and so by (\ref{z11})
\begin{equation}
(1-\gamma)\rho-(3\gamma-1)K\leq\left(\frac{1}{2}-\gamma\right)(g-\lambda)=\left(\dfrac{1}%
{2}-\gamma\right)(\alpha_{s}^{m}+S+K_{s+1}^{m}-\lambda).\label{z22}%
\end{equation}
Putting this into (\ref{z19}) and using (\ref{z6}) we obtain the following
estimate
\begin{eqnarray}
&&(1-\gamma)\alpha_{1}^{r}+(2-3\gamma)a_{r}^{s}+\left(\dfrac{3}{2}-2\gamma\right)\alpha
_{s}^{m}>\frac{1}{2}(K_{s+1}^{m}-\lambda)\label{z23}\\
&&\qquad\qquad +\ (3\gamma-1)K_{1}^{r-1}+\dfrac{5}{2}(2\gamma-1)K_{r}^{m}+(1-\gamma
)h.\nonumber%
\end{eqnarray}
Multiplying  (\ref{z16}) by $(2-3\gamma)>0$ and subtracting from
(\ref{z23}), and noticing that $\frac{5}{2}(2\gamma-1)=2\gamma(2-3\gamma)$ we
obtain
\begin{equation}
(2\gamma-1)(\alpha_{1}^{r}-h)+\left(\gamma-\frac{1}{2}\right)\alpha_{s}^{m}%
+(2-3\gamma)u>\frac{1}{2}(K_{s+1}^{m}-\lambda)+\left(\dfrac{3}{2}-2\gamma
\right)K_{1}^{r-1},\label{z24}%
\end{equation}
and dividing by $2\gamma-1>0$ and using the, equivalent to (\ref{eq}),
equations $(3\gamma+1)(2\gamma-1)=\dfrac{3}{2}-2\gamma$, $(3\gamma+\dfrac
{1}{2})(2\gamma-1)=2-3\gamma$ and $(6\gamma+4)(2\gamma-1)=1$ we obtain the
following basic estimate for the (two tails of the) core measure~$\sigma$:
\begin{equation}
\lbrack a_{1}^{r}-h+(3\gamma+1)(u-K_{1}^{r-1})]+\frac{1}{2}[\alpha_{s}%
^{m}-u+(6\gamma+4)(\lambda-K_{s+1}^{m})]>0,\enspace\ \label{ce1}%
\end{equation}
where we have added and subtracted the term $\frac{1}{2}u$ for reasons that
will become clear in the next section.

This estimate will lead to a contradiction and thus will prove Theorem 1. We
will do this in the following section.

\section{End of the proof of Theorem 1}

Here we will show that both terms in brackets in (\ref{ce1}) must be
nonpositive. This contradicts (\ref{ce1}) and will thus prove Theorem 1.

Consider any measure $\tau$ of the form
\begin{equation}
\tau=\sum_{i=1}^{m}\bar{k}_{i}\delta_{z_{i}},
\end{equation}
where $\bar{k}_{1},\ldots , \bar{k}_{m}>0$ and the $z_{1}<z_{2}<\cdots <z_{m}$ satisfy
the separability inequalities $z_{i+1}-z_{i}>\bar{k}_{i+1}+\bar{k}_{i}$ for
all $1\leq i\leq m-1$ and suppose that the number of positions $m$ in $\tau$
is at most $n$ (the $n$ we have defined in \S 3). The set $E(\tau)$ is
not assumed connected. Consider the set
\begin{equation}
G(\tau)=E(\tau)\backslash\bigcup_{i=1}^{m}I_{i,i}(\tau)\subseteq\lbrack
z_{1},z_{m}]
\end{equation}
that is covered by the nonempty of the intervals $I_{i,j}(\tau)$ where $1\leq
i<j\leq m$. Define the $\bar{K}_{i}^{j}$ similarly to (\ref{a1a}). Then we
have the following.

\specialnumber{13}\proclaim{Lemma}
For every $h$ such that $1<h\leq m$ we have
\begin{equation}
\left|  G(\tau)\cap\lbrack z_{1},z_{h}]\right|  \leq(2\gamma+1)\bar{K}%
_{1}^{h-1}.\label{w1}%
\end{equation}
\endproclaim

\demo{Proof}
The set $G(\tau)\cap\lbrack z_{1},z_{h}]$ is covered by certain intervals
$I_{i,j}(\tau)$ where $1\leq i<j\leq m$. However we know that $I_{i,j}%
(\tau)\subseteq(z_{i},z_{j})$ if $i<j$ and so it would be disjoint from
$[z_{1},z_{h}]$ unless $i<h$. Therefore 
\begin{equation}
G(\tau)\cap\lbrack z_{1},z_{h}]\subseteq\bigcup_{1\leq i<j\leq h-1}%
I_{i,j}(\tau)\cup\bigcup_{1\leq i\leq h-1<j}I_{i,j}(\tau).\label{w2}%
\end{equation}
Consider the measure $\tau^{\prime}=\sum_{i=1}^{h-1}\bar{k}_{i}\delta_{z_{i}}%
$. Then $\bigcup_{1\leq i<j\leq h-1}I_{i,j}(\tau)=G(\tau^{\prime}) $ and since
$\tau^{\prime}$ contains less than $n$ positions we have
\begin{equation}
\left|
{\displaystyle\bigcup_{1\leq i<j\leq h-1}}
I_{i,j}(\tau)\right|  =\left|  G(\tau^{\prime})\right|  \leq2\gamma\bar{K}%
_{1}^{h-1}.\label{w3}%
\end{equation}
Now for the other part consider any interval of the form $I_{i,j}(\tau)$ where
$1\leq i\leq h-1<j$. We have, since $\tau$ satisfies the separability
inequalities,
\begin{equation}
\frak{l}(I_{i,j}(\tau))=z_{j}-\bar{K}_{i}^{j}>z_{h}+\bar{k}_{h}+2\bar{K}%
_{h+1}^{j-1}+\bar{k}_{j}-\bar{K}_{i}^{j}\geq z_{h}-\bar{K}_{1}^{h-1}\hskip.25in
\end{equation}
if $j>h$ and
\begin{equation}
\frak{l}(I_{i,h}(\tau))=z_{h}-\bar{k}_{h}-\bar{K}_{i}^{h-1}\geq z_{h}-\bar
{k}_{h}-\bar{K}_{1}^{h-1}.%
\end{equation}
Therefore since $G(\tau)\cap\lbrack z_{1},z_{h}]\subseteq\lbrack z_{1}%
,z_{h}-\bar{k}_{h}]$ we have
\begin{equation}
G(\tau)\cap\lbrack z_{1},z_{h}]\cap\bigcup_{1\leq i\leq h-1<j}I_{i,j}%
(\tau)\subseteq\lbrack z_{h}-\bar{k}_{h}-\bar{K}_{1}^{h-1},z_{h}-\bar{k}%
_{h}]\hskip.5in\label{w4}%
\end{equation}
and so its measure is at most $\bar{K}_{1}^{h-1}$. Combining (\ref{w4}) with
(\ref{w2}) and (\ref{w3}) we get (\ref{w1}).
\enddemo

{\it Remarks}.  (i) A analogous symmetrical statement holds for $G(\tau
)\cap\lbrack z_{h},z_{m}]$ if $1\leq h<m$.
\vglue4pt
(ii) After Theorem 1 is proved, the above lemma holds for any measure, without
the restriction on the number of positions, and as it can be easily seen is
best possible.
\vglue8pt

Now we can show that both terms in (\ref{ce1}) are nonpositive.

\specialnumber{14}\proclaim{Lemma}
For the core measure $\sigma$ we have
\begin{equation}
\alpha_{1}^{r}-h+(3\gamma+1)(u-K_{1}^{r-1})\leq0.\label{ce2}%
\end{equation}
\endproclaim

\demo{Proof}
We may assume that $r>1$ otherwise there is nothing to prove.

We have by (\ref{z15})
\begin{equation}
u\leq\sum_{i=1}^{r-1}\max(a_{i}+k_{i}-a_{r},0).\label{ce3}%
\end{equation}
Let
\begin{equation}
q=\min\{i:1\leq i\leq r\hbox{ and }a_{i}+k_{i}\geq a_{r}\}.
\end{equation}
(Note that if $q=r$ then $u=0$.) Then using (\ref{ce3}) it is easy to see
that
\begin{equation}
a_{r}-a_{q}+u\leq K_{q}^{r-1}.%
\end{equation}
Therefore we have
\begin{equation}
\alpha_{1}^{r}-h+(3\gamma+1)(u-K_{1}^{r-1})\leq\alpha_{1}^{q}-h-(3\gamma
+1)K_{1}^{q-1}.\hskip.5in\label{ce4}%
\end{equation}
But then, from the considerations in Section 7 and since the definition of $q
$ implies that $I_{i,m+1}(\mu)=\emptyset$ for all $1\leq i<q$, it follows
that the space in $[y_{1},y_{q}]$ not covered by $E(\sigma)$ has measure at
most $h$. Therefore using Lemma 13 we have $\alpha_{1}^{q}-h\leq
(2\gamma+1)K_{1}^{q-1}$, which in view of (\ref{ce4}) easily implies
\pagebreak (\ref{ce2}).
\enddemo

In a similar symmetrical manner we prove that
\begin{equation}
\alpha_{s}^{m}-u+(6\gamma+4)(\lambda-K_{s+1}^{m})\leq0,\label{ce5}%
\end{equation}
noticing that the part of $[y_{s},y_{m}]$ not covered by $E(\sigma)$ has
measure at most $u$ (in view of (\ref{z9})) and using (\ref{z12}).

But now the inequalities (\ref{ce2}) and (\ref{ce5}) contradict the basic core
estimate (\ref{ce1}). Therefore this completes the proof of Theorem 1.

\section{Proof of Theorem 2}

\bigskip It is clearly sufficient to fix a finite positive Borel measure
$\sigma$ and prove (\ref{mi2}) for $\lambda=1$. The functions $F^{+}%
(x)=\sigma((-\infty,x])$ and $F^{-}(x)=\sigma((-\infty,x))$ are measurable as
nondecreasing. Hence for each $h>0$ the set
\begin{equation}
A(h)=\{x:\sigma([x-h,x+h])>2h\}
\end{equation}
is measurable. Letting $E=\{x:M\sigma(x)>1\}$ it is easy to see that
\begin{equation}
E=\bigcup_{h>0}A(h)=\bigcup\{A(h):h\in{\Bbb  Q}\hbox{ and }h>0\}.%
\end{equation}
Hence setting
\begin{equation}
E_{n}=\bigcup\left\{A(h):h\in{\Bbb  Q}\hbox{ and }h>\dfrac{1}{n}\right\}
\end{equation}
we conclude that $E$ is the union of the increasing sequence $(E_{n})$ of
measurable sets. Thus it is enough to show that for any fixed large $n>1$ and
every compact set $K\subseteq E_{n}$ we have
\begin{equation}
\left|  K\right|  \leq C(1+\dfrac{1}{n})\left\|  \sigma\right\| \label{t1}%
\end{equation}
where $C$ is the constant given in (\ref{i4a}).

Fixing $n$ and $K$ as above we can find an interval $[a,b]$ containing $K$ and
such that $b-\sup K,\inf K-a>\left\|  \sigma\right\|  $ and $\sigma
(\{a,b\})=0$ and a partition
\begin{equation}
a=c_{0}<c_{1}<\cdots <c_{N}=b
\end{equation}
of this interval such that
\begin{equation}
\max\limits_{1\leq j\leq N}(c_{j}-c_{j-1})<\dfrac{1}{n^{2}}\hbox{ and }%
\sigma(\{c_{0},c_{1},\ldots , c_{N}\})=0.%
\end{equation}
This is possible since there are at most countably many $x\in{\Bbb  R}$ such
that $\sigma(\{x\})>0$.

Consider now the following positive linear combination of dirac deltas
\begin{equation}
\mu=\sum_{j=1}^{N}\sigma([c_{j-1},c_{j}])\delta_{\frac{c_{j-1}+c_{j}}{2}%
}\; \lower6pt\hbox{.}%
\end{equation}
Then for every $x\in K$ there exists an $h>\dfrac{1}{n}$ such that
\begin{equation}
\sigma([x-h,x+h])>2h.%
\end{equation}
Clearly $h<\left\|  \sigma\right\|  $ and so $[x-h,x+h]\subseteq(a,b)$. Choose
$j$ and $s$ such that $c_{j}<x-h\leq c_{j+1}$ and $c_{s-1}\leq x+h<c_{s}$ and
let $h^{\prime}=\max(c_{s}-x,x-c_{j})>h$. Clearly $h^{\prime}-h<\dfrac
{1}{n^{2}}$. We have
\begin{equation}
\mu([x-h^{\prime},x+h^{\prime}])\geq\sigma([c_{j},c_{s}]\geq\sigma
([x-h,x+h])>2h>2\dfrac{n}{n+1}h^{\prime}\hskip.35in%
\end{equation}
and so
\begin{equation}
K\subseteq\{x:M\mu(x)>\dfrac{n}{n+1}\}
\end{equation}
and so since $\left\|  \mu\right\|  \leq\left\|  \sigma\right\|  $ by applying
Theorem 1 we get (\ref{t1}). This completes the proof of Theorem 2.

\section{Proof of Theorem 3}

To prove Theorem 3 we assume (in view of Theorem 1) that there exists an
admissible positive linear combination of dirac deltas $\mu$ such that
\begin{equation}
\left|  E(\mu)\right|  =2C\left\|  \mu\right\|
\end{equation}
and such that $\left|  E(\nu)\right|  <2C\left\|  \nu\right\|  $ for every
positive linear combination of dirac deltas $\nu$ that contains less positions
than $\mu$, where $C$ is the constant given in (\ref{i4a}).

Now we fix an integer $n>1$ and consider the set
\begin{equation}
\Omega=\{{\bf b}=(y_{1},\ldots , y_{n};k_{1},\ldots , k_{n})\in{\Bbb  R}^{2n}%
:y_{1}\leq\cdots\leq y_{n}\hbox{ and }k_{1},\ldots , k_{n}\geq0\}.%
\end{equation}
Then to every ${\bf b}=(y_{1},\ldots , y_{n};k_{1},\ldots , k_{n})\in\Omega$ we
associate the measure
\begin{equation}
\sigma({\bf b})=\sum_{i=1}^{n}k_{i}\delta_{y_{i}}%
\end{equation}
and the intervals
\begin{equation}
I_{i,j}({\bf b})=[y_{j}-k_{i}-\cdots-k_{j},y_{i}+k_{i}+..+k_{j}]
\end{equation}
for all $1\leq i\leq j\leq n$ (where as usual $[a,b]=\emptyset$ if $b<a$).

Of course the mapping ${\bf b}\rightarrow\sigma({\bf b})$ is not
one-to-one. But it is easy to see (for example using a limiting argument) that
for any measure $\tau=\sum_{i=1}^{m}h_{i}\delta_{z_{i}}$ where $z_{1}%
<\cdots <z_{m}$ and $h_{1},\ldots , h_{m}>0$ and for any ${\bf b}\in\Omega$ such
that $\tau=\sigma({\bf b})$ we have $E(\tau)=\bigcup_{1\leq i\leq j\leq
n}I_{i,j}({\bf b})$.
\vglue4pt
We will use the following well-known lemma.

\specialnumber{15}\proclaim{Lemma}
Let ${\cal C}$ be a finite collection of closed intervals in ${\Bbb  R}$
such that their union $\bigcup{\cal C}$ is an interval $[x,y]$ where $x<y
$. Then there is a subcollection ${\cal C}_{0}=\{[a_{1},b_{1}%
],\ldots , [a_{N},b_{N}]\}$ of ${\cal C}$ such that $\bigcup{\cal C}%
_{0}=\bigcup{\cal C}${\rm ,} satisfying the following
\begin{equation}
a=a_{1}<a_{2}<\cdots <a_{N}=y\label{k1}%
\end{equation}
and
\begin{equation}
a_{2}\leq b_{1}<b_{2},\ldots , a_{N}\leq b_{N-1}<b_{N}.\label{k2}%
\end{equation}
\endproclaim

As it is well known to prove the above lemma it suffices to pick
${\cal C}_{0}$ of minimal cardinality among all subcollections
${\cal C}^{\prime}$ of ${\cal C}$ satisfying $\bigcup{\cal C}%
^{\prime}=\bigcup{\cal C}$, and so no element of ${\cal C}_{0}$ is
contained in any union of other elements of ${\cal C}_{0}$. 
The intervals of ${\cal C}_{0}$ can be arranged so that (\ref{k1}) is satisfied; then
(\ref{k2}) follows easily from the fact that $\bigcup{\cal C}_{0}$ is the interval
$[x,y]$.

Then we will apply the following proposition.

\specialnumber{10}\proclaim{Proposition}
Let $\tau$ be an admissible positive linear combination of dirac deltas
containing exactly $n>1$ positions such that $R(\nu)<R(\tau)$ for every
positive linear combination of dirac deltas $\nu$ that contains less than $n $
positions. Then there exists an admissible measure
\begin{equation}
\tau^{\ast}=\sum_{i=1}^{n}k_{i}^{\ast}\delta_{y_{i}^{\ast}},
\end{equation}
where all $k_{1}^{\ast},\ldots , k_{n}^{\ast}>0$ and all $y_{1}^{\ast}%
<\cdots <y_{n}^{\ast}$ are rational numbers and such that
\begin{equation}
R(\tau^{\ast})\geq R(\tau).%
\end{equation}
\endproclaim

\demo{Proof}
Suppose that $\tau=\sigma({\bf b}_{0})$ where
 $${\bf b}_{0}%
=(y_{1}^{(0)},\ldots , y_{n}^{(0)};k_{1}^{(0)},\ldots , k_{n}^{(0)})\in\Omega$$ is
uniquely determined. By scaling we may assume that
\begin{equation}
E(\tau)=[y_{1}^{(0)}-k_{1}^{(0)},y_{n}^{(0)}+k_{n}^{(0)}]=[0,1].%
\end{equation}
Note that then $k_{1}^{(0)}+\cdots + k_{n}^{(0)}=\left\|  \tau\right\|  \leq1$ for
otherwise $R(\tau)<R(\delta_{0})$.

Now applying Lemma 15 to the collection ${\cal C}=\{I_{i,j}(\tau):1\leq
i\leq j\leq n$ and $I_{i,j}(\tau)\neq\emptyset\}$ we can find a
subcollection $\{I_{i_{1},j_{1}}(\tau),\ldots , I_{i_{N},j_{N}}(\tau)\}$ of
${\cal C}$ that still covers $[0,1]$ and satisfies
\begin{equation}
0=\frak{l}(I_{i_{1},j_{1}}(\tau))<\cdots <\frak{l}(I_{i_{N},j_{N}}(\tau))=1
\end{equation}
and
\begin{equation}
\frak{l}(I_{i_{p},j_{p}}(\tau))\leq\frak{r}(I_{i_{p-1},j_{p-1}}(\tau
))<\frak{r}(I_{i_{p},j_{p}}(\tau))
\end{equation}
for all $p=2,\ldots , N$. It is easy to see that we must have $(i_{1},j_{1})=(1,1)
$ and $(i_{N},j_{N})=(n,n)$. Fixing the set of pairs $\{(i_{1},j_{1}%
),\ldots , (i_{N},j_{N})\}$ we now consider the following set
\begin{eqnarray}
\hskip.5in\Omega^{\ast}&\hskip-5pt=\hskip-5pt&\{ {\bf b}=(y_{1},\ldots ,
y_{n};k_{1},\ldots , k_{n})\in\Omega :y_{1}-k_{1}=0,y_{n}+k_{n}=1,\\
&\hskip-5pt\hskip-5pt&y_{j_{p}}-K_{i_{p}}^{j_{p}}\leq
y_{j_{p+1}}-K_{i_{p+1}}^{j_{p+1}}\hbox{ for all }1\leq p\leq N-1,\nonumber\\
&\hskip-5pt\hskip-5pt&y_{j_{p}}-K_{i_{p}}^{j_{p}}\leq
y_{i_{p-1}}+K_{i_{p-1}}^{j_{p-1}}\leq y_{i_{p}}+K_{i_{p}}^{j_{p}}\hbox{ for all }2\leq
p\leq N\hbox{ and}\nonumber\\ &\hskip-5pt\hskip-5pt&k_{1}+\cdots +
k_{n}\leq1\}.\nonumber
\end{eqnarray}
It is easy to see that $\Omega^{\ast}$ is a nonempty (since ${\bf b}_{0}%
\in\Omega^{\ast}$) compact convex polyhedron contained in a codimension $2$
affine subspace of ${\Bbb  R}^{2n}$. Moreover, it is easy to find nonzero
vectors ${\bf v}_{1},\ldots , {\bf v}_{M}$ such that all the conditions that
define $\Omega^{\ast}$ (including the conditions defining $\Omega$) can be
written as
\begin{equation}
{\bf v}_{1}.{\bf b}=0,\quad {\bf v}_{2}.{\bf b}=1,\quad {\bf v}%
_{3}.{\bf b}\leq0,\ldots , {\bf v}_{M}.{\bf b}\leq0,\label{dr}%
\end{equation}
and, moreover, ${\bf v}_{1}={\bf e}_{1}-{\bf e}_{n+1}$, ${\bf v}%
_{2}={\bf e}_{n}+{\bf e}_{2n}$, and all the entries in all ${\bf v}%
_{1},\ldots , {\bf v}_{M}$ are from the set $\{-1,0,1\}$. Considering the linear
functional $F$ with
\begin{equation}
F({\bf b})=k_{1}+\cdots + k_{n}=({\bf e}_{n+1}+\cdots + {\bf e}_{2n}%
).{\bf b}%
\end{equation}
and applying the standard result from the theory of linear programming we
conclude that there exists an extreme point (vertex) ${\bf b}^{\ast
}=\{y_{1}^{\ast},\ldots , y_{n}^{\ast};k_{1}^{\ast},\ldots , k_{n}^{\ast}\}$ of
$\Omega^{\ast}$ such that
\begin{equation}
F({\bf b}^{\ast})=\min\{F({\bf b):b}\in\Omega^{\ast}\}{\bf \leq}%
F({\bf b}_{0}{\bf )}.%
\end{equation}
Let now $m_{1}=1<m_{2}=2<\cdots <m_{s}\leq M$ be all the indices such that
equality holds in the corresponding relation from (\ref{dr}) when ${\bf b}
$ is replaced by ${\bf b}^{\ast}$. Then it is clear that since
${\bf b}^{\ast}$ is a vertex of $\Omega^{\ast}$ the linear system
\begin{equation}
{\bf v}_{m_{1}}.{\bf b}=0,{\bf v}_{m_{2}}.{\bf b}=1,{\bf v}%
_{m_{3}}.{\bf b}=0,\ldots , {\bf v}_{m_{s}}.{\bf b}=0
\end{equation}
must have ${\bf b}^{\ast}$ as its unique solution and since all
coefficients are integers we conclude that all the $2n$ coordinates of
${\bf b}^{\ast}$ must be {\it rational numbers}.

Consider now the measure $\tau^{\ast}=\sigma({\bf b}^{\ast})$. Since
${\bf b}^{\ast}\in\Omega^{\ast}$ it is easy to see that
\begin{equation}
\lbrack0,1]\subseteq\bigcup_{1\leq i\leq j\leq n}I_{i,j}({\bf b}^{\ast
})=E(\tau^{\ast});
\end{equation}
moreover,
\begin{equation}
\left\|  \tau^{\ast}\right\|  =F({\bf b}^{\ast})\leq F({\bf b}%
_{0})=\left\|  \tau\right\|  .%
\end{equation}
Hence $R(\tau^{\ast})\geq R(\tau)$ and the assumptions on $\tau$ combined with
Lemma 1 (and its proof given in \cite{Mel}) now imply that $\tau^{\ast}$ must
contain exactly $n$ positions and may assumed admissible (without changing its
basic property that all its positions and masses are rational). This completes
the proof of the proposition.
\enddemo

Using now the above proposition we can find an admissible measure $\mu^{\ast}$
whose masses and positions are rational numbers and such that $R(\mu
^{\ast})\geq2C$. But then $R(\mu^{\ast})>2C$ violates Theorem 1 and also
$R(\mu^{\ast})=2C$ leads to a contradiction since $R(\mu^{\ast}) $ must be a
rational number whereas $C$ is irrational. This completes the proof of Theorem 3.

\section{Appendix}

Here we will briefly sketch the construction from \cite{Mel} that leads to the
lower bound in (\ref{i2a}) thus showing that the inequality in Theorem 1 is
actually best possible.

For any admissible measute $\mu$ as in (\ref{a1}) we consider the following
modified norm
\begin{equation}
\left\|  \mu\right\|  ^{\ast}=k_{0}+2k_{1}+\cdots + 2k_{n}+k_{n+1} 
\end{equation}
and the corresponding modified ratio
\begin{equation}
\hbox{\quad}R^{\ast}(\mu)=\frac{\left|  E(\mu)\right|  -k_{0}-k_{n+1}%
}{\left\|  \mu\right\|  ^{\ast}}=\frac{y_{n+1}-y_{0}}{k_{0}+2k_{1}%
+\cdots + 2k_{n}+k_{n+1}}.\hskip.5in
\end{equation}
It is easy to see that $R^{\ast}(\mu)>R(\mu)>1$ for any admissible $\mu$.
Moreover by applying a reflection-translation procedure one can show (see
\cite{Mel}) that for any admissible measure $\mu$ and every $\varepsilon>0$
there exists a measure $\widetilde{\mu}$ such that $R(\widetilde{\mu})\geq
R^{\ast}(\mu)-\varepsilon$. This measures $\tilde{\mu}$ will consist of a
large number of translated copies of $\mu$ (and its symmetric one). Hence any
admissible measure $\mu$ also satisfies $R^{\ast}(\mu)\leq C$.

Then we consider any measure $\nu$ that satisfies the separability condition
(\ref{a2}). We do not assume that $E(\nu)$ is connected. Writing
$\nu$ as $\sum_{i=1}^{n}k_{i}\delta_{y_{i}}$ where $k_{i}>0$ and
$y_{1}<\cdots <y_{n}$, we fix integers $1\leq s,r\leq n$ and define the measure
\begin{equation}
T_{s,r}\nu=k_{0}\delta_{y_{0}}+\nu+k_{n+1}y_{n+1}, %
\end{equation}
where\
\begin{equation}
y_{0}=2y_{1}-y_{s}-2k_{1}-k_{s}\,\hbox{, }\,k_{0}=y_{s}-y_{1}-K_{2}^{s-1}%
\end{equation}
\noindent and
\begin{equation}
y_{n+1}=2y_{n}-y_{r}+2k_{n}+k_{r}\hbox{, }k_{n+1}=y_{n}-y_{r}-K_{r+1}%
^{n-1}.%
\end{equation}

It is easy to show (see \cite{Mel}) that $E(T_{s,r}\nu)$ does not have more gaps
than $E(\nu)$. That is, the added intervals $(y_{0},y_{1}-k_{1}),(y_{n}%
+k_{n},y_{n+1})$ are contained in $E(T_{s,r}\nu)$. Hence the operation
$T_{s,r}$ does not create any new gaps. However we have the advantage of using
the special interval $I_{0,n+1}(T_{s,r}\nu)$, which will be nonempty if $s>t$,
to possibly cover gaps of our initial set $E(\nu)$. For this purpose we argue
as follows.

Let $\mu$ be any, admissible now, measure written as $\mu=\sum_{i=1}^{m}%
k_{i}^{\prime}\delta_{z_{i}}$ where $k_{i}^{\prime}>0$ and $z_{1}<\cdots <z_{n}$
where for simplicity we assume that $z_{1}=0$. Fixing now two positive real
numbers $A,\alpha>0$ we consider the scaled measure $\alpha.\mu$ defined by
\begin{equation}
\alpha.\mu=\sum_{i=1}^{m}\alpha k_{i}^{\prime}\delta_{\alpha z_{i}}.%
\end{equation}
Clearly the measure $\alpha.\mu$ is also admissible and so the measure
\begin{equation}
\nu=\mu+\hbox{trasl}_{A}(\alpha.\mu)=\sum_{i=1}^{m}k_{i}^{\prime}\delta
_{z_{i}}+\sum_{i=1}^{m}\alpha k_{i}^{\prime}\delta_{\alpha z_{i}+A}=\sum
_{i=1}^{n}k_{i}\delta_{y_{i}},\hskip.5in
\end{equation}
where $n=2m$, satisfies the separability inequalities as long as
$A>k_{m}^{\prime}+\alpha k_{1}^{\prime}$. We will next take as $s$ the last
position of $\mu$, so $s=n=2m$, and as $r$ the first position of the
translated $\alpha.\mu$, so $t=m+1$ and consider the measure
\begin{equation}
T\mu=T_{2m,m+1}\nu.\label{tm}%
\end{equation}

Then in \cite{Mel} it is shown that by choosing
\begin{equation}
\alpha=2R(\mu)\hbox{ and }A=(\alpha^{2}-\alpha)\left\|  \mu\right\|
+(\alpha-1)k_{1}^{\prime}%
\end{equation}
the measure $T\mu$ will be admissible (hence $E(T\mu)$ is connected) and
moreover
\begin{equation}
R^{\ast}(T\mu)=\frac{20R(\mu)^{2}-4R(\mu)}{12R(\mu)^{2}-2R(\mu)+1}.%
\end{equation}
\vglue4pt

Let now $f(x)=\dfrac{20x^{2}-4x}{12x^{2}-2x+1}$. Starting from the admissible
measure $\mu_{0}=\delta_{0}+\delta_{3}$ we define the sequence of positive
linear combinations of dirac deltas $(\mu_{p})_{p\geq0}$ (all whose masses and
positions are rational numbers) as follows. Having defined $\mu_{p}$ consider
$T\mu_{p}$ and apply the reflection-translation procedure to obtain a measure
$\mu_{p+1}$ such that $R(\mu_{p+1})\geq R^{\ast}(T\mu_{p})-\varepsilon
_{p}=f(R(\mu_{p}))-\varepsilon_{p}$ where the $\varepsilon_{p}>0$ tend to $0$
sufficiently fast. Then we will have $R(\mu_{p})\rightarrow\dfrac{11+\sqrt
{61}}{12}=1.5675208\ldots$ as $p\rightarrow\infty$. This implies the lower bound
in (\ref{i2a}). After the first few steps these measures will be rather complicated.

However each such measure $\mu_{p}$ will contain a large number of translated
copies of $T\mu_{p-1}$ (and its symmetric one) so it will have a specific
structure. To study this structure let us consider the gap interval
$[a_{0},a_{n+1}]$ of the measure $T\mu$ defined in (\ref{tm}). It is easy to
see that it starts with $a_{1}-a_{0}=\break(1+\alpha)\left\|  \mu\right\|  $
followed by a copy of $J(\mu)$, then by a gap of length $\alpha\left|
J(\mu)\right|  $ (that is completely covered by $I_{0,n+1}(T\mu)$), then by a
copy of $J(\alpha.\mu)=\alpha J(\mu)$ and then by $a_{n+1}-a_{n}%
=\alpha\left\|  \mu\right\|  $. These easily imply that the pair $(J_{0}%
^{+},J_{n+1}^{-})$ has the same structure as the good pairs described in
Section 5 and that both of its intervals are {\it clean}. Moreover its core
is equal to a copy of the measure $a.\mu$ and the $\mu$ corresponds to the
intermediate measure $\nu$ considered in Section~7. Also (assuming all
positions and masses integers), it is easy to see that $\left|  T(J_{0}%
^{+},J_{n+1}^{-})\right|  =\left|  J(T\mu)\right|  $ and $H(J_{0}^{+}%
,J_{n+1}^{-})=\left\|  \mu\right\|  ^{\ast}$; thus their ratio is equal to
$R^{\ast}(T\mu)-1$. So compared with the considerations in Section 7 we
conclude that $T\mu$ shows, in a sense, the tightest possible structure.

In our proof of Theorem 1 we have actually shown that certain measures $\tau$
with $R(\tau)>C$ must have (or can be used to produce) segments that behave in
a structurally similar fashion as the $T\mu$'s. However to prove the sharp
upper bound we had to consider the effect of the more general operator
$T_{s,r}$ with $r<s$ which  makes it necessary to also study certain aspects of
the internal structure of the core, which leads to the basic core estimate
(\ref{ce1}). The fact that in a sense $r$ must be as small as possible and $s$
as large as possible is reflected by the inability to satisfy (\ref{ce1}).
This is what actually leads to the proof of the upper bound.

\end{document}